\documentclass[12pt,reqno]{amsproc}%
\usepackage{amsfonts}
\usepackage{amsmath}
\usepackage{amssymb}
\usepackage{graphicx}%
 \theoremstyle{plain}

\newtheorem{corollary}{Corollary}[section]

\newtheorem{definition}{Definition}[section]
\newtheorem{example}{Example}[section]

\newtheorem{lemma}{Lemma}[section]

\newtheorem{proposition}{Proposition}[section]
\newtheorem{remark}{Remark}[section]

\newtheorem{theorem}{Theorem}[section]
 \numberwithin{equation}{section}
\begin{document}

\centerline{\Large\bf{Some Examples of Blackadar and Kirchberg's MF Algebras}}

\vspace{1cm}

\centerline{Don Hadwin \qquad \qquad
 and \qquad \qquad  Junhao
Shen\footnote{The  author is partially supported by an NSF grant.}}
\bigskip

\centerline{\small{Department of Mathematics and Statistics,
University of New Hampshire, Durham, NH, 03824}}

\vspace{0.2cm}

\centerline{Email: don@math.unh.edu  \qquad \quad  and\qquad \quad
\    jog2@cisunix.unh.edu \qquad }

\bigskip

\noindent\textbf{Abstract: } In the paper,  we provide some examples
of Blackadar and Kirchberg's   MF algebras by considering  minimal
and maximal tensor products of MF algebras and crossed products of
MF algebras by finite groups or an integer  group. We also present
some examples of C$^*$-algebras, whose BDF extension semigroups are
not   groups. These examples include, for example, $
 C_r^*(F_n)\otimes_{max} C^*(F_n)  $, $ C_r^*(F_n)\otimes_{min}
C^*(F_n) $, $C_r^*(H_1\ast H_2)$ with $2\le |H_1|<\infty$ and $3\le
|H_2|<\infty$ where $|H_1|$, $|H_2|$ are the orders of the groups
$H_1$, or $H_2$ respectively, and several others.

\vspace{0.2cm} \noindent{\bf Keywords:} MF algebras, tensor products, crossed product

\vspace{0.2cm} \noindent{\bf 2000 Mathematics Subject
Classification:} Primary 46L10, Secondary 46L54

 \section{Introduction}


 The notion of MF algebra was introduced   by Blackadar and Kirchberg in \cite{BK}.    A separable C$^*$-algebra
  is an MF algebra if and only if it can be embedded into
   $\prod_{k}\mathcal M_{n_k}(\Bbb C)/\sum_{k}\mathcal M_{n_k}(\Bbb C)$ for a sequence of positive
   integers $n_k$, $k=1,2,\ldots$.  It follows   from the definition that a  separable
   quasidiagonal C$^*$-algebra is an MF algebra. 


The class of MF algebras is closely connected to Brown, Douglas and
Fillmore's extension theory (see \cite{BDF}). In \cite{Haag},
  Haagerup and Thorbj{\o}rnsen solved a long standing open problem by showing that $Ext(C^*_r(F_n))$ is not a
  group for $n\ge 2$, where   $C^*_r(F_n)$
   is the reduced C$^*$-algebra of the free group $F_n$ and $Ext(C^*_r(F_n))$ is the BDF's extension
   semigroup. Their result follows from a combination of Voiculescu's argument in \cite{V4} and their remarkable work on proving
   that $C^*_r(F_n)$ is an MF algebra.  Basing on Haagerup and Thorbj{\o}rnsen's approach to BDF's extension semigroup,
   one would be interested in
   finding other natural examples (besides the ones in \cite{An}, \cite{Haag} and \cite{Wa}) of C$^*$-algebras whose BDF extension semigroups are
   not
   groups. This is one of  the motivations of our investigation on MF
   algebras. In \cite{HaSh4}, we considered the full free products
   of MF algebras and showed that $Ext(C_r^*(F_n)\ast_{\Bbb
   C}\mathcal A)$ is not a group when $\mathcal A$ is an MF algebra
   and $n\ge 2$ is a  positive integer. In this paper, we  consider  tensor products of MF algebras
    and crossed products of MF algebras by finite groups or an integer group. In fact,  we are able to show that, for example,
   both $ Ext(C_r^*(F_n)\otimes_{max} C^*(F_n)) $ and $Ext(C_r^*(F_n)\otimes_{min} C^*(F_n))$ are not groups
    (see Corollary 3.2, Corollary 3.3, Corollary 4.1 and Example 3.1).

Blackadar and Kirchberg's MF algebra is also closely connected to
Voiculescu's topological free entropy dimension. In \cite{Voi}, for
a family of self-adjoint elements $x_1,\ldots, x_n$ in a unital
C$^*$-algebra $\mathcal A$, Voiculescu introduced the notion of
topological free entropy dimension of $x_1,\ldots,x_n$.  In the
definition of
  topological free entropy dimension, it requires that
Voiculescu's norm microstate space of $x_1,\ldots,x_n$ is
``eventually" nonempty, which is equivalent to say that the
C$^*$-subalgebra generated by $x_1,\ldots, x_n$ in $\mathcal A$
should be an MF algebra. Thus it is important  to determine which
C$^*$-algebra belongs to the class of MF algebra, where Voiculescu's
topological free entropy dimension is well defined  (More discussion
on topological free entropy dimension can be found in \cite{HaSh2},
\cite{HaSh3}, \cite{HaSh4}).

More connections between MF algebra and the classification of
C$^*$-algebras can also be found in \cite {BK}, \cite{BK2},
\cite{BK3}.

Now we  outline   the main results obtained in the paper. First,
we study the minimal and maximal tensor products of MF algebras
and obtain the following results:

\vspace{0.2cm}
 \noindent{\bf Theorem 3.1: } Suppose $C^*(F_n)$ is the full C$^*$-algebra of the free group $F_n$ ($n\ge 2$) and  $\mathcal B$ is a unital separable MF algebra. Then
$C^*(F_n)\otimes_{max} \mathcal B$ is an MF algebra.

\vspace{0.2cm}

\vspace{0.2cm}
 \noindent{\bf Proposition 3.1 \& 3.2: } Suppose  $\mathcal A$ and  $\mathcal B$ are
  separable MF algebras. If $\mathcal A$ is either exact or quasidiagonal, then $ \mathcal A\otimes_{min} \mathcal B$ is an MF algebra.

\vspace{0.2cm}

We also consider the crossed product of an MF algebra by finite
groups or an integer group and obtain an analogue of Pimsner and
Voiculescu's result \cite{PiVoi} in the context of MF algebra.

\vspace{0.2cm}
 \noindent{\bf Theorem 4.2: }
Suppose that $\mathcal A$ is a finitely generated unital MF algebra  and
$\alpha$ is a homomorphism from $\Bbb Z$ into $Aut(\mathcal A)$ such
that there is a sequence of integers $0\le n_1<n_2<\cdots$
satisfying $ \lim_{j\rightarrow \infty}\|\alpha(n_j) a-a\|=0 $  for
any $a\in \mathcal A$. Then $\mathcal A\rtimes_{\alpha}\Bbb Z $ is
an MF algebra.

\vspace{0.2cm}

Combining with Haagerup and Thorbj{\o}rnsen's result that
$C_r^*(F_n)$ is an MF algebra, we   obtain more examples of
C$^*$-algebras whose BDF extension semigroups are not   groups.

\vspace{0.2cm}
 \noindent{\bf Corollary 3.2: }
Suppose $n,m\ge 2$ are positive integers and $\mathcal B$ is a
unital separable MF algebra. Then $(C_r^*(F_n)\ast_{\Bbb C}\mathcal
B)\otimes_{max}C^*(F_m)$ is an MF algebra and
$Ext((C_r^*(F_n)\ast_{\Bbb C}\mathcal B)\otimes_{max}C^*(F_m))$ is
not a group.

\vspace{0.2cm}
 \noindent{\bf Corollary 3.3: }
Suppose that $\mathcal B  $ is  a  unital separable MF algebras.
Then $Ext(C^*_r(F_n)\otimes_{min}\mathcal B)$ is not a group for
$n\ge 2$.

\vspace{0.2cm}
 \noindent{\bf Corollary   4.2: }
Assume $n, m$ are positive integers and $G$ is a
  finite group with $|G|\ge 2$, where $|G|$ is the order of group $G$. Then   $Ext(C_r(F_n\ast (G\times F_m)))$ is
not a group.

\vspace{0.2cm}
 \noindent{\bf Corollary   4.3: }
Assume that $n\ge 2$ is a positive integer and $H_1, \ldots, H_n$ is
a family of finite groups. Let
$$
H= H_1\ast H_2 \ast \cdots \ast H_n.
$$  Moreover, if there are $1\le i\ne j\le n$ such that
$$ |H_i|\ge 2 \qquad and \qquad |H_j|\ge 3,$$ where $|H|$ denotes the order
of the group $H$, then $Ext(C_r^*(H))$ is not a group.

\vspace{0.2cm}
 \noindent{\bf Corollary 4.4: }
Let $SL_2(\Bbb Z)$ to be the spacial linear group of $2\times 2$
matrices with integer entries. Then $C_r^*(SL_2(\Bbb Z))$ is an MF
algebra and $Ext(C_r^*(SL_2(\Bbb Z)))$ is not a group.

\vspace{0.2cm}

 \noindent{\bf Corollary 4.5:}
Let $C_r^*(F_2)$ be the reduced C$^*$-algebra of the free group
$F_2$.  Let $u_1,u_2$ be   canonical unitary generators of
$C_r^*(F_2)$ and $0<\theta  <1$ be a positive number. Let $\alpha$
be a homomorphism from $\Bbb Z$ into $Aut(C_r^*(F_2))$ induced by
the following mapping: $\forall \ n\in \Bbb Z$,
$$
\alpha(n)(u_1)=e^{2n \pi    \theta \cdot  i} u_1 \qquad and \qquad
\alpha(n)(u_2)=e^{2n\pi   \theta \cdot i} u_2.
$$ Then $C_r^*(F_2)\rtimes_{\alpha}\Bbb Z $ is an MF algebra and
$Ext(C_r^*(F_2)\rtimes_{\alpha}\Bbb Z)$ is not a group.

\vspace{0.2cm}

 We hope our
investigation on MF algebra will be helpful in the study of
Blackadar and Kirchberg's MF algebras, BDF's extension semigroup
and Voiculescu's topological free entropy theory.

The organization of the paper is as follows. In section 2,  we
introduce some notation and preliminaries. In section 3, we study
  tensor products of MF algebras. In
section 4, we consider  crossed product of an MF algebra by an
action of a finite group or an integer group.

 \section{Notation and Preliminaries}

 Suppose $\mathcal H$ is a separable complex Hilbert space and $B(\mathcal H)$ is the set of all
 bounded linear operators on $\mathcal H$.  Suppose $ \{x, x_k\}_{k=1}^\infty$ is a family of elements in $B(\mathcal H)$.
 We say $x_k\rightarrow x$ in $*$-SOT ($*$-strong operator topology) if and only if $x_k\rightarrow x$ in SOT and
 $x_k^*\rightarrow x^*$ in  SOT. 

 For a  sequence of  C$^*$ algebras $\mathcal A_k\subseteq B(\mathcal H_k)$, $k=1,2,\ldots,$  we introduce the direct product of $\mathcal A_k$ as follows.
 $$
 \prod_k \mathcal A_k  = \left\{ \langle x_k\rangle_{k=1}^\infty \ | \ \sup_k\|x_k\|<\infty, \ \ where \ \ x_k\in\mathcal A_k\right\}, $$ where  the norm of any element $\langle x_k\rangle_{k=1}^\infty$ in $
 \prod_k \mathcal A_k$ is defined by
 $$
 \|\langle x_k\rangle_{k=1}^\infty\|= \sup_k\|x_k\|.
 $$ Let $$
 \sum_k \mathcal A_k  = \left\{ \langle x_k\rangle_{k=1}^\infty\in  \prod_k \mathcal A_k \ |
 \ \limsup_{k\rightarrow\infty}\|x_k\|=0, \ \ where \ \ x_k\in\mathcal A_k\right\}.
 $$
 It is easy to see that $\sum_k \mathcal A_k$ is a closed two-sided idea of the   C$^*$-algebra $
 \prod_k \mathcal A_k$. Thus $\prod_k \mathcal A_k/\sum_k \mathcal A_k$ is also a   C$^*$-algebra.
  Denote by $\pi$ the canonical quotient mapping from $\prod_k \mathcal A_k$ onto $\prod_k \mathcal A_k/\sum_k \mathcal A_k$, and denote
 $\pi(\langle x_k\rangle_{k=1}^\infty)$ by $[\langle x_k\rangle_{k=1}^\infty] $ for every element $\langle x_k\rangle_{k=1}^\infty$ in $
 \prod_k \mathcal A_k$.

Suppose $\mathcal A$ is a    C$^*$-algebra on a Hilbert space
$\mathcal H$. Let $\mathcal H^\infty=\oplus_{\Bbb N} \mathcal H$,
   and for every $x\in \mathcal A$, let
$x^\infty$ be the element $\langle x,x,x,\ldots
\rangle=\oplus_{\Bbb N} x $ in $\prod_k \mathcal A^{(k)}\subseteq
B(\mathcal H^\infty)$, where $\mathcal A^{(k)}$ is the $k$-th
copy of $\mathcal A$. Let $\mathcal A^\infty=\{x^\infty \in
\prod_k \mathcal A^{(k)}\ | \ x\in\mathcal A\}$ be a
C$^*$-subalgebra of $\prod_k \mathcal A^{(k)}$. 

 For each positive integer $n$,   denote the $n\times n$ matrix algebra over the complex numbers by
 $\mathcal M_n(\Bbb C)$ and the group of unitary matrices in $\mathcal M_n(\Bbb C)$ by
 $\mathcal U_n(\Bbb C)$.

Recall the definition of residually finite dimensional
C$^*$-algebras as follows.
  \begin{definition} A separable C$^*$-algebra $\mathcal A$ is residually finite dimensional if there is an embedding from $\mathcal A$ into
$\prod_k\mathcal M_{n_k}(\Bbb C)$ for some positive integers $\{n_k\}_{k=1}^\infty$.
\end{definition}

Recall the definition of quasidiagonal C$^*$-algebras as follows.
  \begin{definition} A set of elements $\{a_1,\ldots,a_n\}\subseteq B(\mathcal H)$ is quasidiagonal
if
there  is an increasing sequence of
finite-rank projections $\{p_i\}_{i=1}^\infty$ on $H$ tending
strongly to the identity such that $\| a_j p_i-p_i a_j \|\rightarrow 0$
as $i\rightarrow \infty$ for any $1\le j\le n.$ A  separable C$^*$-algebra
 $\mathcal A\subseteq B(\mathcal H)$ is
quasidiagonal    if there  is an increasing sequence of
finite-rank projections $\{p_i\}_{i=1}^\infty$ on $H$ tending
strongly to the identity such that $\| x p_i-p_i x \|\rightarrow 0$
as $i\rightarrow \infty$ for any $ x \in   \mathcal A.$     
 An abstract separable
C$^*$-algebra $\mathcal A$ is quasidiagonal if there is a faithful
$*$-representation  $\pi: \mathcal A\rightarrow B(\mathcal H)$ on
a Hilbert space $\mathcal H$ such that $\pi(\mathcal A)\subseteq
B(\mathcal H)$ is quasidiagonal.  \end{definition}
 Let   $\Bbb C\langle
X_1,\ldots, X_n\rangle $ be the set of all  noncommutative
polynomials in the indeterminates $X_1,\ldots$, $X_n,X_1^*,\ldots,
X_n^*$. The following lemma is well known.
\begin{lemma}
Suppose  $\mathcal A\subseteq B(\mathcal H)$ is a unital separable
 quasidiagonal C$^*$-algebra  and $x_1,\ldots, x_n$ are elements
in  $\mathcal A$. For any $\epsilon>0$, a finite subset
$\{P_1,\ldots, P_r\}$ of $\Bbb C\langle X_1,\ldots, X_n\rangle$, and
a finite subset $\{\xi_1\ldots, \xi_s\}$ of $\mathcal H$, there is a
finite rank projection $p $ in $B(\mathcal H)$ such that:
\begin{enumerate}
\item [(i)] $\| P_i(p \ x_1p , \ldots, p \ x_np )\xi_j- P_i( x_1 , \ldots,
 x_n )\xi_j\|<\epsilon;$  for all $1\le i\le r$ and $1\le   j\le s.$
\item [(ii)] $|\|P_i(p \ x_1p , \ldots, p \ x_np )\|_{B(p\mathcal H)}-\|P_i( x_1 , \ldots,
 x_n )\|_{\mathcal A} |<\epsilon,$ for all $1\le i\le r.$ \end{enumerate}
\end{lemma}
The following well known result of Voiculescu   plays a crucial role in the study of quasidiagonal C$^*$-algebras.
\begin{lemma}
Suppose $\mathcal A$ is a separable quasidiagonl C$^*$-algebra
and $\pi$ is an essentially faithful $*$-representation of
$\mathcal A$ on a separable Hilbert space $\mathcal H$, i.e.
$\pi(\mathcal A)\cap \mathcal K=\{0\}$ where $\mathcal K$ is the
set of compact operators on $\mathcal A$. Then $\pi(\mathcal
A)\subseteq B(\mathcal H)$ is  quasidiagonal.
\end{lemma}



 The notion of MF algebra was introduced by Blackadar and Kirchberg in \cite{BK}.
  \begin{definition}
  A separable  C$^*$-algebra $\mathcal A$ is an MF algebra if there is an embedding from $\mathcal A$ into $\prod_k\mathcal M_{n_k}(\Bbb C)/\sum_k\mathcal M_{n_k}(\Bbb C)$ for some positive integers $\{n_k\}_{k=1}^\infty.$ \end{definition}

\begin{remark} The following implication is obvious from the definitions: If a C$^*$-algebra $\mathcal A$ is  residually finite dimensional then it is a quasidiagonal C$^*$-algebra. If $\mathcal A$ is a quasidiagonal C$^*$-algebra, then it is an MF algebra.
\end{remark}

 The following lemma gives   equivalent definitions of an MF algebra.

\begin{lemma}
Suppose that $\mathcal A$ is a    unital C$^*$-algebra generated
by a family of   elements $x_1,\ldots, x_n$ in $\mathcal A$. Then
the following are equivalent:
\begin{enumerate}
   \item $\mathcal A$ is an MF algebra;
   \item For any $\epsilon>0$ and any finite family $\{P_1,\ldots,P_r\}$ in $\Bbb C\langle X_1,\ldots,X_n\rangle$, there are a   positive integers $k$ and family of  matrices
   $\{A_1 ,\ldots, A_n \}$ in $\mathcal M_{ k} (\Bbb C)$, such that
   $$
   \max_{1\le j\le r}\left |\|P_j(A_1^{(k)},\ldots, A_n^{(k)})\|- \|P_j(x_1,\ldots,x_n)\|\right |<\epsilon.
   $$
   \item Suppose $\pi: \mathcal A \rightarrow B(\mathcal H)$ is a faithful $*$-representation
    of $\mathcal A$ on an infinite dimensional separable complex Hilbert space $\mathcal H$.
    Then there are a sequence of positive integers $\{m_k\}_{k=1}^\infty$, families of self-adjoint  matrices
   $\{A_1^{(k)},\ldots, A_n^{(k)}\}$ in $\mathcal M_{m_k} (\Bbb C)$ for $k=1,2,\ldots$, and unitary operators $U_k: \mathcal
   H \rightarrow  (\Bbb C^{m_k})^\infty$ for $k=1,2,\ldots$,  such that
  \begin{enumerate} \item $
   \lim_{k \rightarrow \infty}\|P(A_1^{(k)},\ldots, A_n^{(k)})\|= \|P(x_1,\ldots,x_n)\|,   \ \forall \ P\in \Bbb C\langle X_1,\ldots,X_n\rangle;
   $ \item $
      U_k^*\ (A_i^{(k)})^\infty  \ U_k  \rightarrow \pi(x_i) \ in \ *-SOT  \ as \ k \rightarrow \infty,\ \ for \ 1\le i\le n,
   $ \\  where $(A_i^{(k)})^\infty=  A_i^{(k)}\oplus A_i^{(k)}\oplus A_i^{(k)}  \cdots \in B((\Bbb C^n)^\infty)$. \end{enumerate}

   \item Suppose $\pi: \mathcal A \rightarrow B(\mathcal H)$ is a faithful $*$-representation of $\mathcal A$ on an infinite dimensional separable complex Hilbert space $\mathcal H$. Then there is a family of elements $\{a_1^{(k)},\ldots, a_n^{(k)}\}_{k=1}^\infty\subseteq B(\mathcal H)$ such that
        \begin{enumerate}
        \item For each $k\ge 1$, $\{a_1^{(k)},\ldots, a_n^{(k)}\}\subseteq B(\mathcal H)$ is quasidiagonal;
        \item $
   \lim_{k \rightarrow \infty}\|P(a_1^{(k)},\ldots, a_n^{(k)})\|= \|P(x_1,\ldots,x_n)\|,   \ \forall \ P\in \Bbb C\langle X_1,\ldots,X_n\rangle;
   $ \item $ \forall \ 1\le i\le n, \ \
       a_i^{(k)}  \rightarrow  \pi(x_i),  \ in \ *-SOT  \ as \ k \rightarrow \infty.$
        \end{enumerate}
\end{enumerate}

\end{lemma}
\begin{proof}
$(1)\Leftrightarrow (2) \Leftrightarrow (3)$ was proved in
\cite{HaSh4}. $(3)\Rightarrow (4)$ is immediate. We only need to
show $(4)\Rightarrow (2)$. For any $\epsilon>0$ and
$\{P_1,\ldots,P_r\}$ in $\Bbb C\langle X_1,\ldots,X_n\rangle$, by
condition (b), there is some $m\ge 1$ such that
$$
   \max_{1\le j\le r}\left |\|P_j(a_1^{(m)},\ldots, a_n^{(m)})\|- \|P_j(x_1,\ldots,x_n)\|\right
   |<\epsilon/2;
   $$
 and $a_1^{(m)},\ldots, a_n^{(m)}\subseteq B(\mathcal H)$ is quasidiagonal. Thus from Lemma 2.1 it follows that there is a finite rank projection $p$ such that $$
   \max_{1\le j\le r}\left |\|P_j(pa_1^{(m)}p,\ldots, pa_n^{(m)}p)\|_{B(p\mathcal H)}- \|P_j(a_1^{(m)},\ldots, a_n^{(m)})\|\right |<\epsilon/2.
   $$ Hence
  $$
   \max_{1\le j\le r}\left |\|P_j(pa_1^{(m)}p,\ldots, pa_n^{(m)}p)\|_{B(p\mathcal H)}- \|P_j(x_1,\ldots,x_n)\|\right |
   <\epsilon.
   $$  Let $k=dim(B(p\mathcal H))$. Note that $B(p\mathcal H)\simeq \mathcal M_k(\Bbb C)$. Let
   $A_i=pa_i^{(m)}p$ for $1\le i\le n$. Then we have the desired result.
\end{proof}

The following result   can be found in \cite{Br} (See also
\cite{VoiQusi},\cite{Haag}) which connects
  MF algebras to BDF's extension semigroups.

\begin{lemma}   If $\mathcal A$ is a unital separable  MF algebra and $Ext(\mathcal A)$ is a group, then $\mathcal A$ is quasidiagonal.

\end{lemma}
\begin{proof}
For the purpose of completeness, we sketch its proof here. By
Voiculescu's characterization of quasidiagonality, to show that
$\mathcal A$ is quasidiagonal, it suffices to show that for any
$a_1,\ldots, a_n$ in $\mathcal A$ and any $\epsilon>0$, there is a
unital completely positive map $\psi:\mathcal A\rightarrow \mathcal
M_k(\Bbb C)$ for some positive integer $k$ such that (1)
$\|\psi(a_i)\|\ge \|a\|-\epsilon$; and (2)
$\|\psi(a_ia_j)-\psi(a_i)\psi(a_j)\|\le \epsilon.$

Note that $\mathcal A$ is an MF algebra. There is a  unital
embedding $\rho: \mathcal A \rightarrow \prod_k \mathcal
M_{n_k}(\Bbb C)/\sum_k \mathcal M_{n_k}(\Bbb C)$ for some positive
integers $\{n_k\}_{k=1}^\infty$. Let $\mathcal H = \Bbb C^{n_1}
\oplus \Bbb C^{n_2}\oplus \cdots$ be a Hilbert space and $p_k$ the
projection from $\mathcal H$ onto the subspace $0\oplus\cdots \oplus
0 \oplus \Bbb C^{n_k}\oplus 0 \oplus \cdots$. Let $\mathcal
K(\mathcal H)$ be the set of all compact operators on $\mathcal H$
and $\pi$   the quotient map from $B(\mathcal H)$ onto $B(\mathcal
H)/\mathcal K(\mathcal H)$. Then
$$
\rho: \mathcal A \rightarrow \prod_k \mathcal M_{n_k}(\Bbb C)/\sum_k
\mathcal M_{n_k}(\Bbb C) \hookrightarrow B(\mathcal H)/\mathcal
K(\mathcal H).
$$
Let $[\langle A_i^{(k)}\rangle_{k=1}^\infty]=\rho(a_i)\in  \prod_k
\mathcal M_{n_k}(\Bbb C)/\sum_k \mathcal M_{n_k}(\Bbb C) $, and
$x_i= \langle A_i^{(k)}\rangle_{k=1}^\infty \in \prod_k \mathcal
M_{n_k}(\Bbb C)  \hookrightarrow B(\mathcal H)$ for $1\le i\le n$.
Thus  (i) $p_kx_i=p_kx_ip_k$.  By Proposition 3.3.5 in \cite{BK},
we can assume that  (ii) $\|a_i\|= \lim_k \|p_kx_ip_k\|=\lim_k
\|A_i^{(k)}\|  $; and (iii) $ \|a_ia_j\|=\lim_k
\|p_kx_ip_kx_jp_k\|=\lim_k \|A_i^{(k)}A_j^{(k)} \|$ for $1\le
i,j\le n.$

Since $Ext(\mathcal A)$ is a group, by \cite{Ar}, there is a unital
completely positive map $\phi:\mathcal A \rightarrow B(\mathcal H)$
such that $\rho= \pi\circ \phi$. It follows that
$\pi(\phi(a_i))=\pi(x_i)$ and $\pi(\phi(a_ia_j))=\pi(x_ix_j)$ for
$1\le i,j\le n$. Then $\phi(a_i )-x_i$ and $ \phi(a_ia_j)-x_ix_j$
are in $\mathcal K(\mathcal H)$, whence (iv) $\lim_k\|p_kx_ip_k
-p_k\phi(a_i)p_k\|=0$ and (v) $\lim_k\|p_kx_ix_jp_k
-p_k\phi(a_ia_j)p_k\|=0$

Define $$\psi_k:\mathcal A \rightarrow B(p_k\mathcal H)\simeq
\mathcal M_{n_k}(\Bbb C)\qquad  by \qquad  \psi_k(a)= p_k\phi(a)p_k,
\ \  \forall \ a\in\mathcal A.
$$ Then each $\psi_k$ is a unital completely positive map from
$\mathcal A$ to $\mathcal M_{n_k}(\Bbb C)$. When $k$ is large
enough, by  (ii) and (iv), we know that $\psi_k$ satisfies (1).
By (i), (iv) and (v),  $\psi_k$ satisfies (2). Hence $\mathcal A$
is quasidiagonal.
\end{proof}
We also need the following result by Rosenberg \cite{Ro}.
\begin{lemma}
Suppose that $G$ is a countable discrete group. If $C^*_r(G)$ is quasidiagonal, then $G$ is an amenable group.

\end{lemma}

\section{Tensor products of MF Algebras}

 Suppose that $\mathcal A$ and $\mathcal B$ are two C$^*$-algebras. Denote the minimal tensor product of $\mathcal A$ and $\mathcal B$ by $\mathcal A\otimes_{min}\mathcal B$ and the maximal tensor product of $\mathcal A$ and $\mathcal B$ by $\mathcal A\otimes_{max}\mathcal B$.

 Suppose that $F_n$ is the free group on $n$ many generators with $n\in \Bbb N $.
 Let $C^*_r(F_n)$ be the reduced C$^*$-algebra associated with $F_n$ and  $C^*(F_n)$
 be the full C$^*$-algebra of the group $F_n$.

 \subsection{Maximal Tensor Product } In this subsection, we consider the question whether the maximal
  tensor product  of $C^*(F_n)$ and an  MF algebra is again an MF algebra. Before we present our main result, we will introduce a few lemmas first.

 \begin{lemma}
 For any $\delta>0$, there is a polynomial $P\in C\langle X\rangle$ with $P(0)=0$  such that,
    $$\sup_{0\le t\le 1}\|\sqrt t-P(t)\|\le \delta.$$
 \end{lemma}
\begin{proof}
It is an easy exercise.
\end{proof}

\begin{lemma}
Suppose that $\mathcal H$ is a  Hilbert space. For any $\delta>0$,
there is a   function  $$f_\delta : (0,\infty) \rightarrow
(0,\infty)$$ such that the following hold:
\begin{enumerate}
  \item $\lim_{t\rightarrow 0} f_\delta(t)\le \delta;$ 
  \item Suppose   $p$ and $q$ are  mutually orthogonal, equivalent, finite rank,  projections in $B(\mathcal H)$
   with a partial isometry $w\in B(\mathcal H)$ satisfying
      $w^*w=p, ww^*=q$.   
      Suppose    $ u_1,\ldots, u_n,\\ v_1,\ldots v_m $ are   unitary elements in $B(\mathcal H)$
      such that, for some $t>0$,
      $$
   u_iv_j=v_iu_j,\qquad    and \qquad  \|pv_j-v_jp\|\le
   t, \quad \forall \ 1\le i\le n, 1\le j\le m.
      $$ Let
 $$
\begin{aligned}
U_i&=pu_ip+ \sqrt{p-pu_ipu_i^*p} \ w^*+ w \
\sqrt{p-pu_i^*pu_i p}-wpu_i^*pw^* ,\quad \forall \  1\le i\le n\\
V_j&=pv_jp+wv_jw^*,\qquad \forall \ 1\le j\le m.
\end{aligned}
 $$   Then
   \begin{enumerate}
   \item $U_1,\ldots, U_n$ are uniatry elements in $B((p+q)\mathcal
   H)$;
 \item     $\forall \ 1\le i\le n$,
 $1\le j\le m$, $$
    \|U_iV_j- V_jU_i\|\le f_\delta(t).
      $$
      \end{enumerate}
\end{enumerate}

\end{lemma}
\begin{proof}
 It is trivial to verify that
$U_1,\ldots, U_n$ are unitary elements in $B((p+q)\mathcal H)$. 


Note $u_iv_j=v_ju_i$ and $\|pv_j-v_jp\|\le t$. We know that
\begin{equation}
\|pv_jp\cdot pu_i^*pu_ip-pu_i^*pu_ip\cdot pv_jp\|\le 3t .
\end{equation} Similarly,
\begin{equation}
\|pv_jp\cdot pu_ipu_i^*p-pu_ipu_i^*p\cdot pv_jp\|\le 3t
\end{equation}

 For such $\delta>0$, by Lemma 3.1, there is a polynomial $P\in\Bbb
C\langle X\rangle $ such that
$$ \max_{0\le \theta\le 1}|\sqrt{\theta}-P(\theta)|\le \delta/4.  $$
This implies
$$\|\sqrt{p-pu_i^*pu_ip} - P(p-pu_i^*pu_ip)\| \le \delta/4, \qquad \|\sqrt{p-pu_i pu_i^*p} - P(p-pu_i^*pu_ip)\| \le  \delta/4.$$
For such polynomial $P$, by (3.1) and (3.2),  there is some
positive number $D_\delta>0$ (only depending on $P$, not on $t$)
so that
\begin{align}
\|pv_jpP(p-pu_i^*pu_ip)-P(p-pu_i^*pu_ip) pv_jp\|\le D_\delta \cdot t
, \notag
\end{align} and
\begin{align}\|pv_jpP(p-pu_ipu_i^*p)-P(p-pu_ipu_i^*p) pv_jp\|\le D_\delta \cdot  t
.\notag
\end{align}
Then
\begin{align}
\|pv_jp\sqrt{p-pu_i^*pu_ip}-\sqrt{p-pu_i^*pu_ip}\  pv_jp\|\le
D_\delta\cdot t +\delta /2  ,
\end{align} and
\begin{align}\|pv_jp\sqrt{p-pu_i pu_i^*p}-\sqrt{p-pu_i pu_i^*p}\  pv_jp\|\le D_\delta \cdot t +\delta/2.
\end{align}
 Since
$$
V_j =pv_jp+wv_jw^*,\qquad \forall \ 1\le j\le m,
$$ from (3.3) and (3.4) we have \begin{align}
\|U_iV_j-V_jU_i\|&= \| pu_ipv_jp+ \sqrt{p-pu_ipu_i^*p}\  v_j  \ w^*+ w \ \sqrt{p-pu_i^*pu_ip}\ pv_jp-wpu_i^* pv_j w^*  \notag\\
&\quad  \ \ \ \ -(pv_jpu_ip+ pv_jp\sqrt{p-pu_ipu_i^*p} \ w^*+ w v_j  \ \sqrt{p-pu_i^*pu_ip}-w v_jp u_i^*pw^*)\|\notag\\
&\le 4t+ 2tD_\delta + \delta.\end{align} Let $f_\delta(t)=4t+
2tD_\delta + \delta$. By the construction of $f_\delta$, we know
that $f_\delta$ only depends on $\delta$, not on $u_i,v_j,p,q$ and
$t$. Thus we have desired result.
\end{proof}

\begin{lemma}
Suppose that $\mathcal A$ is a  unital MF algebra generated by a
family of unitary elements $x_1,\ldots,x_n$, and $\mathcal B$ is a
unital C$^*$ algebra generated by $y_1,\ldots, y_m$.  For any
$\epsilon>0$, any
 $\{P_1,\ldots, P_r\}\subseteq \Bbb C\langle X_1,\ldots, X_n,
Y_1,\ldots,Y_n\rangle$  and $ \{Q_1,\ldots, Q_t\}\subseteq \Bbb
C\langle X_1,\ldots,X_n\rangle, $  there is a quasidiagonal
C$^*$-algebra $\mathcal A_1$  with a family of unitary generators
$\{z_1,\ldots, z_n\}$ 
 so
that the following hold:
\begin{enumerate}
   \item For any $1\le j\le r$, $$ \begin{aligned}
  \|P_j(x_1\otimes 1,\ldots, &x_n\otimes 1, 1\otimes y_1,\ldots, 1\otimes y_m)\|_{\mathcal A\otimes_{max}\mathcal B}\\
  & \le  \|P_j(z_1\otimes 1,\ldots,z_n\otimes 1, 1\otimes  y_1,\ldots,
   1\otimes y_m)\|_{\mathcal A_1\otimes_{max}\mathcal B}.\end{aligned}
   $$ \item For any $1\le j\le t$,
   $$
 \left | \|Q_j(  z_1 ,\ldots,    z_n)\|_{\mathcal A_1} -  \|Q_j(    x_1,\ldots,
  x_n)\|_{\mathcal A} \right | \le \epsilon.$$
\end{enumerate}
\end{lemma}

\begin{proof}
Note that $\mathcal A$ is an MF algebra and $x_1,\ldots, x_n$ are
unitary elements in $\mathcal A$. Thus there are a sequence of
positive integers $\{t_l\}_{l=1}^\infty$ and families of unitary
matrices
   $\{B_1^{(l)},\ldots, B_n^{(l)}\}$ in $\mathcal M_{t_l} (\Bbb C)$ for $l=1,2,\ldots$, such that
   $$
   \lim_{l \rightarrow \infty}\|Q(B_1^{(l)},\ldots, B_n^{(l)})\|= \|Q(x_1,\ldots,x_n)\|,   \ \forall \ Q\in \Bbb C\langle X_1,\ldots,X_n\rangle.$$


 For any positive integers $N_1 $ and
$1\le i\le n$, we define
$$
  \begin{aligned}
     D(N_1;i) & = B_i^{({N_1})}\oplus B_i^{({N_1+1})}\oplus B_i^{({N_1+2})}\oplus \cdots\\
  \end{aligned}
$$
It is not hard to see that, for any $Q\in \Bbb C\langle
X_1,\ldots,X_n\rangle$,
 \begin{align}
\|Q(D(N_1;1)  ,\ldots, D(N_1;n)  )\| &=\sup_{k\ge N_1}
\|Q(B_1^{(k)},\ldots, B_n^{(k)})\|
 \ge \|Q(x_1,\ldots,x_n)\|\\
 \lim_{N_1\rightarrow \infty}\|Q(D(N_1;1) ,\cdots, D(N_1;n) )\|& =
 \|Q(x_1,\ldots,x_n)\|
\end{align}
Therefore, for the chosen set $\{Q_1,\ldots Q_t\}$, by (3.7) there
exists some $N_1$ such that
\begin{align}
\left | \|Q_j(D(N_1;1) ,\ldots, D(N_1;n) )\| -
 \|Q_j(x_1,\ldots,x_n)\| \right |\le \epsilon, \ \ for \ \ 1\le j\le t.
\end{align}

Let $z_i= D(N_1;i)$ for $1\le i\le n$ and $\mathcal A_1$ be the
unital C$^*$-algebra generated by $z_1,\ldots,z_n$. Hence, by
(3.8)
\begin{equation}
\left | \|Q_j(z_1  ,\ldots,  z_n  )\| -
 \|Q_j(x_1,\ldots,x_n)\| \right |\le \epsilon, \ \ for \ \ 1\le j\le t.
\end{equation} It is trivial to see that $z_1,\ldots, z_n$ are unitary elements in $\mathcal A_1$.
 Moreover by inequality (3.6), there is a $*$-homomorphism
$$
 \rho_1: \mathcal A_1\rightarrow \mathcal  A,
$$ such that
$\rho_1(z_i)=  x_i  $ for $1\le i\le n.$
 By the
property of the maximal tensor product, it induces a
$*$-homomorphism
$$
\rho_1\otimes_{max} id : \ \mathcal A_1\otimes_{max}\mathcal B
\rightarrow \mathcal A \otimes_{max}\mathcal B.
$$ In follows that, for any
$1\le j\le r$,   \begin{align}
  \|P_j(x_1\otimes 1,\ldots, &x_n\otimes 1, 1\otimes y_1,\ldots, 1\otimes y_m)\|_{\mathcal A\otimes_{max}\mathcal B}\notag\\
  & \le  \|P_j(z_1\otimes 1,\ldots,z_n\otimes 1, 1\otimes  y_1,\ldots,
   1\otimes y_m)\|_{\mathcal A_1\otimes_{max}\mathcal B}.\end{align}
   Combining with   (3.9) and (3.10), we have completed the proof of the
  Lemma.

\end{proof}

\begin{lemma}  Suppose $C^*(F_n)$ is the full C$^*$-algebra
of the free group $F_n$ and $u_1,\ldots, u_n$ are   canonical
generators of $C^*(F_n)$. Suppose $\mathcal B$ is a unital
quasidiagonal C$^*$-algebra generated by unitary elements
$z_1,\ldots,z_m$.   Suppose $\rho$ is a faithful representation of $
C^*(F_n) \otimes_{max} \mathcal B $ on a separable Hilbert space
$\mathcal H$.

For any $\epsilon >0$, any   $\{P_1,\ldots,P_r\} $
of $ \Bbb C\langle X_1,\ldots, X_n,Y_1,\ldots, Y_m\rangle$,
  and any   $\{Q_1,\ldots,Q_t\}$ in $\Bbb
C\langle X_1,\ldots,X_n\rangle$, there are a positive integer
$k_1$ and  $k_1\times k_1$ matrices
$$
\{U_1,\ldots,U_n\}\subseteq \mathcal U_{k_1}(\Bbb C) \qquad
\{V_1,\ldots,V_m\}\subseteq \mathcal M_{k_1}(\Bbb C)
$$ such that (1)
$$
\|U_iV_j-V_jU_i\|\le \epsilon, \qquad \forall \ \ 1\le i\le n, \
1\le j\le m;
$$ (2), $\forall\  1\le j\le t$,
$$
\left| \|Q_j(V_1,\ldots,V_m)\|- \|Q_j(z_1,\ldots,z_m)\|\right|\le
\epsilon;
$$ and (3), $\forall\ 1\le i\le r$,
$$
  \| P_i( u_1\otimes 1,\ldots,u_n\otimes 1,
1\otimes z_1,\ldots,1\otimes z_m)\|_{ C^*(F_n)\otimes \mathcal B}
\le \|P_i(U_1,\ldots,U_n,V_1,\ldots,V_m)\|+\epsilon.
$$
\end{lemma}
\begin{proof}
Since $\rho$ is a faithful representation of
$C^*(F_n)\otimes_{max}\mathcal B$, without loss of generality  we
can further assume that $\rho$ is an essentially faithful
representation. Then  $\rho(1\otimes \mathcal B )\subseteq
B(\mathcal H)$ is quasidiagonal. It follows that there is a
sequence of finite rank projections $ p_k$, $k=1,2,\ldots,$ on
$\mathcal H$ such that (i)
 $1\le j\le m$,  $$  \|p_k \rho(1\otimes z_j) - \rho(1\otimes z_j)p_k\|
\rightarrow 0; $$ and (ii) $ \forall\ 1\le j\le t$,
\begin{align} \|Q_j( p_k \rho( 1\otimes z_1) p_k ,\ldots, p_k
\rho(1 \otimes z_m) p_k )\| \rightarrow  \|Q_j( \rho(1\otimes z_1)
,\ldots, \rho(1 \otimes z_m))\| .
\end{align}
Let $$ t_k =\max_{  1\le j\le m} \{   \|p_k
\rho(1\otimes z_j) - \rho(1\otimes z_j)p_k\|\}.
$$ Then $t_k$ tends to $0$.  Let $q_k$ be a  projection in $B(\mathcal H)$ so that $p_k, q_k$ are mutually
orthogonal, equivalent projections, and let $w_k$ be the partial
isometry so that  $w_k^*w_k=p_k, w_k w_k^*=q_k$. Let
\begin{align}
U_i^{(k)}&=p_k\rho(u_i\otimes 1)p_k+ \sqrt{p_k-p_k\rho(u_i\otimes
1)p_k\rho(u_i\otimes 1)^*p_k} \ w_k^* \notag \\
&\ \ \ + w_k \
\sqrt{p_k-p_k\rho(u_i\otimes 1)^*p_k\rho(u_i\otimes 1)p_k}-w_kp_k\rho(u_i\otimes 1)^*p_kw_k^* ,\quad \forall \  1\le i\le n\\
V_j^{(k)}&=p_k\rho(1\otimes z_j)p_k+w_k\rho(1\otimes
z_j)w_k^*,\qquad \forall \ 1\le j\le m
  \end{align}
   be  elements in
$B((p_k+q_k)\mathcal H)$.  For   $\delta=\frac
\epsilon 4>0 $, let
$$f_\delta: (0,\infty) \rightarrow (0,\infty)$$  be the function as in Lemma 3.2. Then
\begin{enumerate}
\item[(a)] $ \lim_{t\rightarrow 0} f_\delta(t)\le  \epsilon/4 $;  and  \item
[(b)]   $ U_i^{(k)}$ is a  unitary element in $B((p_k+q_k)\mathcal
H)$; and \item [(c)] $\forall \ 1\le i\le n, \ 1\le j\le m$
     \begin{align}
    \|U_i^{(k)} V_j^{(k)}-   V_j^{(k)}U_i^{(k)} \|\le f_\delta(t_k).
    \end{align} \end{enumerate}
We observe that $U_i^{(k)}, V_j^{(k)}$, as $k$ goes to infinity,
converge in strong operator topology to
 $$
 \begin{pmatrix}
   \rho(u_i\otimes 1) & 0\\
   0 &-  \rho(u_i\otimes 1)^*
 \end{pmatrix} \qquad and \qquad \begin{pmatrix}
   \rho(1\otimes z_j) & 0\\
   0 &\rho(1\otimes z_j)
 \end{pmatrix}.
 $$
It follows that, for any $1\le i\le r$,
 \begin{align}
\liminf_{k\rightarrow \infty}&\| P_i(U_1^{(k)},\ldots,   U_n^{(k)},V_1^{(k)},
\ldots, V_m^{(k)})\|\notag\\
 &\ge \| P_i( \rho(u_1\otimes 1),\ldots,   \rho(u_n\otimes 1),
\rho(1\otimes z_1),\ldots,\rho(1\otimes z_m))   \| \notag\\
&= \|P_i( u_1\otimes 1,\ldots,u_n\otimes 1,  1\otimes
z_1,\ldots,1\otimes z_m) \|_{  C^*(F_n)  \otimes_{max}\mathcal B
}
\end{align}
 On the other hand,
 \begin{align}
\lim_{k} &\left | \|Q_j(V_1^{(k)}, \ldots, V_n^{(k)})\| - \|Q_j(
z_1,\ldots,
  z_m)\|_{\mathcal B}  \right | \notag \\
  & =\lim_{k}
\left | \|Q_j( p_k \rho(1\otimes z_1)p_k ,\ldots,
    p_k \rho(1\otimes z_m) p_k )\|   - \|Q_j(
z_1,\ldots,
  z_m)\|_{\mathcal B}   \right | \notag \\
&=\left | \|Q_j(   \rho(1\otimes z_1)  ,\ldots,
     \rho(1\otimes z_m)  )\|   - \|Q_j(
z_1,\ldots,
  z_m)\|_{\mathcal B}  \right | \tag{By (3.11)}  \\
  & =\left | \|Q_j(   1\otimes z_1   ,\ldots,
      1\otimes z_m   )\|_{C^*(F_n)\otimes_{max}\mathcal B}   - \|Q_j(
z_1,\ldots,
  z_m)\|_{\mathcal B}  \right |   \notag\\
   &=0.
\end{align}
By (3.14), (3.15),  and (3.16), we know when $k$ is large enough,
there are a positive integer $k_1=dim((p_{k_1}+q_{k_1})\mathcal
H)$ and $k_1\times k_1$ matrices
$$
\{U_1,\ldots,U_n\}\subseteq \mathcal U_{k_1}(\Bbb C),\qquad
\{V_1,\ldots,V_m\}\subseteq \mathcal M_{k_1}(\Bbb C)
$$ such that (i)
$$
\|U_iV_j-V_jU_i\|\le  \epsilon, \qquad \forall \ \ 1\le i\le n, \
1\le j\le m;
$$ (ii) $\forall\  1\le j\le t$,
$$
\left| \|Q_j(V_1,\ldots,V_m)\|- \|Q_j(z_1,\ldots,z_m)\|_{\mathcal
B}\right|\le \epsilon;
$$ and (iii), $\forall\ 1\le i\le r$,
$$
  \| P_i( u_1\otimes 1,\ldots,u_n\otimes 1,
1\otimes z_1,\ldots,1\otimes z_m)\|_{ C^*(F_n)
\otimes_{max}\mathcal B  }\le
\|P_i(U_1,\ldots,U_n,V_1,\ldots,V_m)\|+\epsilon.
$$
\end{proof}

Recall that $C^*(F_n)$ is the full C$^*$-algebra of the free
group $F_n$ and $u_1,\ldots, u_n$ is a family of canonical
generators of  $C^*(F_n)$. We assume that $\mathcal B$ is an MF
algebra with a family of unitary generators $v_1,\ldots, v_m$.
From now on till the end of this subsection, we will assume that
$\{P_r\}_{r=1}^\infty$, or $\{Q_t\}_{t=1}^\infty$ respectively,
is the collection of all polynomials in $\Bbb C\langle
X_1,\ldots, X_n,Y_1,\ldots, Y_m\rangle$, or $\Bbb C\langle
Y_1,\ldots, Y_m\rangle$ respectively, with the rational
coefficients.
\begin{lemma} Suppose  $u_1,\ldots, u_n$ is a family of canonical
generators of $C^*(F_n)$. Suppose $\mathcal B$ is an MF algebra with
a family of unitary generators $v_1,\ldots, v_m$. Let
$\{P_r\}_{r=1}^\infty, \{Q_t\}_{t=1}^\infty$ be as above.

 For any $r_0\ge 1$, there is
some $r_1\ge 1$ so that the following hold: for  any positive
integer $k\ge 1$,   if
$$
\{U_1, \ldots, U_n\}\subseteq \mathcal U_k(\Bbb C)\qquad \{
V_1,\ldots, V_m\}\in \mathcal M_k(\Bbb C),
$$ is a family of $k\times k$ matrices satisfying,
$$
\|U_iV_j-V_jU_i\|\le \frac 1 {r_1 } \ , \qquad \ \   \forall \ 1\le
i\le n, \ 1\le j\le m ;
$$ and $$
\left |\|Q_j(V_1,\ldots, V_m)\|-\|Q_j(v_1,\ldots, v_m)\|_{\mathcal
B}\right |\le \frac 1 {r_1 } \ , \qquad \forall  \ 1\le j\le r_1,
$$ then, $ \forall \  1\le i\le r_0$
$$
 \|P_i(U_1,\ldots,U_n,V_1,\ldots,V_m) \| \le
\|P_i(u_1\otimes 1 ,\ldots,u_n\otimes 1, 1\otimes
v_1,\ldots,\otimes v_m)\|_{C^*(F_n)\otimes_{max}\mathcal B}+\frac
1 {r_0}.
$$

\end{lemma}
\begin{proof}
We will prove the lemma by using contradiction. Suppose the result
of the lemma does not hold. There are some positive integer
$r_0>0$, a sequence of positive integers $\{k_r\}_{r=1}^\infty,$
and sequences of matrices
$$
\{U_1^{(r)}, \ldots, U_n^{(r)}\}\subseteq \mathcal U_{k_r}(\Bbb
C)\qquad \{ V_1^{(r)},\ldots, V_m^{(r)} \}\in \mathcal
M_{k_r}(\Bbb C),
$$ satisfying \begin{enumerate}
\item [(i)] $$ \|U_i^{(r)}V_j^{(r)}-V_j^{(r)}U_i^{(r)}\|\le \frac 1 { r } \ , \qquad \ \   \forall \
1\le i\le n, \ 1\le j\le m ;
$$ \item [(ii)] and $$
\left |\|Q_j(V_1^{(r)},\ldots, V_m^{(r)})\|-\|Q_j(v_1,\ldots,
v_m)\|_{\mathcal B}\right |\le \frac 1 { r } \ , \qquad \forall \
1\le j\le r;
$$ \item [(iii)] but
\begin{align}
\max_{1\le i\le r_0}
\{\|P_i(U_1^{(r)},\ldots,U_n^{(r)},V_1^{(r)},\ldots,V_m^{(r)}) \|
- \|P_i(u_1\otimes 1 ,\ldots,u_n\otimes 1, 1\otimes
v_1,\ldots,\otimes v_m)\|_{C^*(F_n)\otimes_{max}\mathcal
B}\}>\frac 1 {r_0}.
\end{align}
\end{enumerate}
Consider the C$^*$-algebra $\prod_r \mathcal M_{k_r}(\Bbb
C)/\sum_r \mathcal M_{k_r}(\Bbb C)$. Let
$$
 \begin{aligned}
    x_i&= [\langle U_i^{(r)}\rangle_{r=1}^\infty],  \qquad \forall \ 1\le i\le n\\
    y_j&= [\langle V_j^{(r)}\rangle_{r=1}^\infty],  \qquad \forall \ 1\le j\le m
 \end{aligned}
$$ be the elements in $\prod_r \mathcal M_{k_r}(\Bbb C)/\sum_r
\mathcal M_{k_r}(\Bbb C)$. It is not hard to check:
\begin{enumerate} \item [(iv)] $x_1,\ldots,x_n$ are unitary element; \item [(v)] $x_iy_j=y_jx_i$ for
all $1\le i\le n,$ $1\le j\le m$; \item [(vi)] for all $j\ge 1$,
$$\| Q_j(y_1,\ldots, y_m)\|_{\prod_r \mathcal M_{k_r}(\Bbb
C)/\sum_r \mathcal M_{k_r}(\Bbb C)}=\|Q_j(v_1,\ldots,
v_m)\|_{\mathcal B}
$$\end{enumerate}
By the universal property of the C$^*$-algebra
$C^*(F_n)\otimes_{max}\mathcal B$, we know there is $*$-homomorphism
$$\psi: C^*(F_n)\otimes_{max}\mathcal B \rightarrow C^*(x_1,\ldots,
x_n,y_1,\ldots,y_m)$$ such that
$$
\psi(u_i\otimes 1)=x_i, \qquad \psi(1\otimes v_j)= y_j, \qquad
\forall \ 1\le i\le n, 1\le j\le m,
$$ where $C^*(x_1,\ldots,
x_n,y_1,\ldots,y_m)$ is the unital C$^*$-algebra generated by
$x_1,\ldots, x_n, y_1,\ldots, y_m$ in $\prod_r \mathcal M_{k_r}(\Bbb
C)/\sum_r \mathcal M_{k_r}(\Bbb C)$. Hence, for all $i\ge 1$,
$$\begin{aligned}
\limsup_r\|P_i(U_1^{(r)},\ldots,U_n^{(r)},&V_1^{(r)},\ldots,V_m^{(r)})
\|  = \|P_i(x_1,\ldots,x_n,y_1,\ldots, y_m)\|\\& \le
\|P_i(u_1\otimes 1 ,\ldots,u_n\otimes 1, 1\otimes v_1,\ldots,\otimes
v_m)\|_{C^*(F_n)\otimes_{max}\mathcal B}.\end{aligned}
$$
This contradicts with the inequality (3.17). The proof of the lemma
is completed.
\end{proof}

Now we are able to show our main result in this section.
\begin{theorem}
Suppose  $n$ is a positive integer and $\mathcal B$ is a finitely
generated unital MF algebra. Then $C^*(F_n)\otimes_{max} \mathcal
B$ is an MF algebra.
\end{theorem}
\begin{proof}
Suppose that $u_1,\ldots, u_n$ are   canonical unitary generators of
$C^*(F_n)$ and $v_1,\ldots, v_m$ are unitary generators of $\mathcal
B$. Let $r_0, t\ge 1$, $ 0<\epsilon<1/r_0$,
$\{P_i\}_{i=1}^{r_0}\subset \Bbb C\langle X_1, \ldots,X_n,
Y_1,\ldots,Y_m\rangle$ and $\{Q_j\}_{j=1}^t \subset  \Bbb C\langle
Y_1,\ldots,Y_m\rangle$. By Lemma 3.3, there is a quasidiagonal
C$^*$-algebra $\mathcal B_1$ with a family of unitary generators
$\{z_1,\ldots, z_n\}$ so that the following hold:
\begin{enumerate}
   \item [(i)] For any $1\le i\le {r_0}$, $$ \begin{aligned}
  \|P_i(u_1\otimes 1,\ldots, &u_n\otimes 1, 1\otimes v_1,\ldots, 1\otimes v_m)\|_{C^*(F_n)\otimes_{max}\mathcal B}\\
  & \le  \|P_i(u_1\otimes 1,\ldots,u_n\otimes 1, 1\otimes  z_1,\ldots,
   1\otimes z_m)\|_{C^*(F_n)\otimes_{max}\mathcal B_1}.\end{aligned}
   $$ \item [(ii)] For any $1\le j\le t$,
   $$
 \left | \|Q_j(  v_1 ,\ldots,    v_m)\|_{\mathcal B} -  \|Q_j(    z_1,\ldots,
  z_m)\|_{\mathcal B_1} \right | \le \epsilon/2.$$
  \end{enumerate}
  Let $\rho$ be a faithful
$*$-representation $\rho: C^*(F_n)\otimes_{max} \mathcal
B_1\rightarrow B(\mathcal H)$ on a separable Hilbert space
$\mathcal H$.  By lemma 3.4, there are a positive integer $k$ and
$k\times k$ matrices
$$
\{U_1,\ldots,U_n\}\subseteq \mathcal U_k(\Bbb C) \qquad
\{V_1,\ldots,V_m\}\subseteq \mathcal M_{k}(\Bbb C)
$$ such that \begin{enumerate} \item [(iii)] \begin{align}
\|U_iV_j-V_jU_i\|\le \epsilon, \qquad \forall \ \ 1\le i\le n, \
1\le j\le m;
\end{align}  \item [(iv)]  $\forall\  1\le j\le t$,
$$
\left| \|Q_j(V_1,\ldots,V_m)\|- \|Q_j(z_1,\ldots,z_m)\|_{\mathcal B_1}\right|\le
\epsilon/2;
$$  \item [(v)]  $\forall\ 1\le i\le {r_0}$,
$$
  \| P_i( u_1\otimes 1,\ldots,u_n\otimes 1,
1\otimes z_1,\ldots,1\otimes z_m)\|_{ C^*(F_n)\otimes_{max}
\mathcal B_1 }\le
\|P_i(U_1,\ldots,U_n,V_1,\ldots,V_m)\|+\epsilon/2.
$$\end{enumerate}

By (ii) and (iv) we know that  $\forall\  1\le j\le t$,
\begin{align}
\left| \|Q_j(V_1,\ldots,V_m)\|- \|Q_j(v_1,\ldots,v_m)\|_{\mathcal
B }\right|\le \epsilon;
\end{align}
By (i) and (v), we know that $\forall\ 1\le i\le {r_0}$,
\begin{align}
  \| P_i( u_1\otimes 1,\ldots,u_n\otimes 1,
1\otimes v_1,\ldots,1\otimes v_m)\|_{ C^*(F_n)\otimes_{max}
\mathcal B}\le \|P_i(U_1,\ldots,U_n,V_1,\ldots,V_m)\|+\epsilon.
\end{align}

 By Lemma
3.5 and inequalities (3.18), (3.19)we know that, when $\epsilon$
is small enough and $t$ is large enough,  for $1\le i\le
r_0$\begin{align}
  \| P_i( u_1\otimes 1,\ldots,u_n\otimes 1,
1\otimes v_1,\ldots,1\otimes v_m)\|_{ C^*(F_n)\otimes_{max}
\mathcal B }\ge \|P_i(U_1,\ldots,U_n,V_1,\ldots,V_m)\|- 1/r_0,
\end{align} Combining with (3.20), we obtain
\begin{align}
 \left | \| P_i( u_1\otimes 1,\ldots,u_n\otimes 1,
1\otimes v_1,\ldots,1\otimes v_m)\|_{ C^*(F_m)\otimes_{max}
\mathcal B }- \|P_i(U_1,\ldots,U_n,V_1,\ldots,V_m)\|\right |\le
1/r_0.
\end{align}

By (3.18) and (3.22), we know that, for any $r_0\in\Bbb N$, there
are  $k$ and matrices
$$
\{U_1,\ldots,U_n\}\subseteq \mathcal U_k(\Bbb C) \qquad
\{V_1,\ldots,V_m\}\subseteq \mathcal M_{k}(\Bbb C)
$$  satisfying
$$
\|U_iV_j-V_jU_i\|\le 1/ r_0, \qquad \forall \ \ 1\le i\le n, \
1\le j\le m;
$$   and   $\forall\ 1\le i\le r_0$,
$$
 \left | \| P_j( u_1\otimes 1,\ldots,u_n\otimes 1,
1\otimes v_1,\ldots,1\otimes v_m)\|_{ C^*(F_m)\otimes_{max} \mathcal
B }- \|P_j(U_1,\ldots,U_n,V_1,\ldots,V_m)\|\right |\le  1/r_0.
$$ By Lemma 2.3, $C^*(F_n)\otimes_{max} \mathcal B$ is an MF
algebra.
\end{proof}
\begin{remark}
The argument in the proof of Theorem 3.1 also applies to the case
when $\mathcal B$ is countably generated. Thus the preceding
result can be generalized as
 follows: Suppose that $n$ be a positive integer and $\mathcal B$ is a  unital separable MF
algebra. Then $C^*(F_n)\otimes_{max} \mathcal B$ is an MF algebra.
\end{remark}

The following result by Brown and Ozawa (see Proposition 7.4.5 of
\cite{BrOz}) follows directly from the preceding theorem.
\begin{corollary}
Suppose $n,m\ge 2$ are positive integers. Then
$C^*(F_n)\otimes_{max}C^*(F_m)$ is an MF algebra, i.e. there is an
embedding $$ \rho: C^*(F_n)\otimes_{max}C^*(F_m) \rightarrow
\prod_{k} \mathcal M_{n_k}(\Bbb C)/\sum_{k} \mathcal M_{n_k}(\Bbb C)
$$ for some positive integers $\{n_k\}_{k=1}^\infty$.
\end{corollary}

Note   Haagerup and Thorbj{\o}rnsen proved that $C_r^*(F_n)$ is an
MF algebra. By Lemma 2.4 and 2.5 and Theorem 3.1 and Corollary 5.3
in \cite{HaSh4}, we have the following result.
\begin{corollary}
Suppose $n,m\ge 2$ are positive integers and $\mathcal B$ is a
separable MF algebra. By Corollary 5.3 in \cite{HaSh4},
$C_r^*(F_n)\ast_{\Bbb C}\mathcal B$ is an MF algbebra, where
$C_r^*(F_n)\ast_{\Bbb C}\mathcal B$ is the unital full free
product of $C_r^*(F_n) $ and $ \mathcal B$. Then
$(C_r^*(F_n)\ast_{\Bbb C}\mathcal B)\otimes_{max}C^*(F_m)$ is an
MF algebra and $Ext((C_r^*(F_n)\ast_{\Bbb C}\mathcal
B)\otimes_{max}C^*(F_m))$ is not a group by Lemma 2.4 and Lemma
2.5.
\end{corollary}


\subsection{Minimal Tensor Product }  The following result might have been known to experts.
For the purpose of completeness, we also include it here.
\begin{proposition}
Suppose that $\mathcal A$ and $\mathcal B$ are two separable MF algebras. If $\mathcal A$ is exact, then $\mathcal A\otimes_{min}\mathcal B$ is also an MF algebra.
\end{proposition}
\begin{proof}
Since $\mathcal B$ is an MF algebra, there is an embedding
$$
\rho: \mathcal B \rightarrow \prod_k\mathcal M_{n_k}(\Bbb C)/\sum_k\mathcal M_{n_k}(\Bbb C)
$$ for a family of positive integers $\{n_k\}_{k=1}^\infty$.  This embedding $\rho$ induces another embedding
$$
id\otimes_{min}\rho : \mathcal A\otimes_{min}\mathcal B \rightarrow \mathcal A\otimes_{min}\left ( \prod_k\mathcal M_{n_k}(\Bbb C)/\sum_k\mathcal M_{n_k}(\Bbb C)\right).
$$ Note $\mathcal A$ is an exact C$^*$-algebra. We have
$$
\mathcal A\otimes_{min} \left (\prod_k\mathcal M_{n_k}(\Bbb C)/\sum_k\mathcal M_{n_k}(\Bbb C)\right)\simeq \left ( \mathcal A\otimes_{min} \prod_k\mathcal M_{n_k}(\Bbb C)\right )/\left (\mathcal A\otimes_{min}\sum_k\mathcal M_{n_k}(\Bbb C)\right).
$$ Since $\mathcal A$ is an MF algebra, by Corollary 3.4.3 in \cite{BK}, any separable C$^*$-subalgebra of
$$\left ( \mathcal A\otimes_{min} \prod_k\mathcal M_{n_k}(\Bbb C)\right )/\left (\mathcal A\otimes_{min}\sum_k\mathcal M_{n_k}(\Bbb C)\right)\subseteq
  \prod_k \mathcal A\otimes_{min}\mathcal M_{n_k}(\Bbb
C) /  \sum_k  \mathcal A\otimes_{min}\mathcal M_{n_k}(\Bbb C)
$$ is an MF algebra. Therefore,  $\mathcal A\otimes_{min}\mathcal B$ is also an MF algebra.
\end{proof}

Note   Haagerup and Thorbj{\o}rnsen proved that $C_r^*(F_n)$ is an
MF algebra.  By Lemma 2.4 and 2.5 and preceding proposition,  we
have the following result.
\begin{corollary}
Suppose that $\mathcal B$ is a unital  separable MF algebra. Then
$C^*_r(F_n)\otimes_{min}\mathcal B$ is an MF algebra for any $n\ge
2$. Moreover, $Ext(C^*_r(F_n)\otimes_{min}\mathcal B)$ is not a
group by Lemma 2.4 and Lemma 2.5.
\end{corollary}

When both two C$^*$ algebras was not exact, we have the following result.
\begin{proposition}
Suppose  $\mathcal A$ is a separable quasidiagonal C$^*$-algebra and   $\mathcal B$ a  separable MF algebras. Then $\mathcal A\otimes_{min}\mathcal B$ is   an MF algebra.
\end{proposition}
\begin{proof}
We   need only to show that if $\mathcal A_1$, and $\mathcal B_1$
respectively, is a finite generated C$^*$-subalgebra of $\mathcal
A$, and $\mathcal B$ respectively, then  $\mathcal
A_1\otimes_{min}\mathcal B_1$ is MF algebra. Suppose that
$x_1,\ldots, x_n,$ and $y_1,\ldots,y_m$ respectively, is a family
of generators
 of $\mathcal A_1$, and $\mathcal B_1$ respectively.
 By Lemma 2.3, to show  $\mathcal A_1\otimes_{min}\mathcal B_1$ is MF algebra, it suffice to show the following: {\em for any $\epsilon>0$ and any finite subset $\{P_1,\ldots,P_r\} $ of noncommutative $*$-polynomials
 $\Bbb C\langle X_1,\ldots,X_n,Y_1,\ldots,Y_m\rangle$, there are a positive integer $k$ and matrices
 $$
 \{A_1,\ldots,A_n,B_1,\ldots,B_m\}\subseteq \mathcal M_{k}(\Bbb C)
 $$ such that
 $$
  \max_{1\le j\le r} \left | \|P_j(A_1,\ldots,A_n,B_1,\ldots,B_m)\|- \|P_j(x_1\otimes 1,\ldots,x_n\otimes 1, 1\otimes y_1,\ldots,1\otimes y_m)\|_{\mathcal A\otimes_{min}\mathcal B}  \right |\le \epsilon.
 $$}

  Let $\rho$ be a   faithful $*$-representation of $\mathcal B$ on a separable Hilbert space $\mathcal H_1$.
 Let $\pi$ be an essentially faithful $*$-representation of $\mathcal A$ on a separable Hilbert space $\mathcal H_2$. By Lemma 2.2, we know that $\pi(\mathcal A)\subseteq B(\mathcal H_2)$ is quasidiagonal.
 Then there is a sequence of finite rank projections $q_s$, $s=1,2,\ldots,$ in $B(\mathcal H_2)$ such that (i) $q_s \rightarrow I$ in SOT and (ii) $\|q_s \pi(x) -\pi(x)q_s\|\rightarrow 0$ for any $x\in \mathcal A$.
  Let $\mathcal H=\mathcal H_2\otimes \mathcal H_1$ be a Hilbert space. By (i), we know that
  $$\begin{aligned}
   \langle q_s \pi(x_1) q_s\otimes 1, &\ldots, q_s \pi(x_n) q_s\otimes 1, q_s\otimes \rho(y_1),\ldots,
   q_s\otimes \rho(y_m) \rangle \rightarrow^{*-SOT} \notag\\
  &\qquad \qquad \qquad \langle  \pi(x_1)  \otimes 1, \ldots,   \pi(x_n)  \otimes 1, 1\otimes \rho(y_1),\ldots, 1\otimes \rho(y_m) \rangle .\end{aligned}
  $$Thus, for $  \ 1\le j\le r$,
  \begin{align}
  \|P_j(&\pi(x_1)\otimes 1, \ldots,\pi(x_n) \otimes 1, 1\otimes\rho( y_1),\ldots,1\otimes
   \rho(y_m))\|\notag\\& \le \liminf_{s\rightarrow \infty} \|P_j(q_s \pi(x_1) q_s\otimes 1,
   \ldots, q_s \pi(x_n) q_s\otimes 1, q_s\otimes \rho(y_1),\ldots, q_s\otimes \rho(y_m) )\|.
   \end{align} By (ii), we know that,  for $  \ 1\le j\le r$,
  \begin{align}
  \lim_{s\rightarrow \infty} & \left \|P_j(q_s \pi(x_1) q_s\otimes 1,  \ldots, q_s \pi(x_n) q_s\otimes 1,
  q_s\otimes \rho(y_1),\ldots,
   q_s\otimes \rho(y_m)) \right .
 \notag \\
  &  \left .  - (q_s\otimes 1) \left (P_j(\pi(x_1)\otimes 1,\ldots,\pi(x_n) \otimes 1, 1\otimes\rho( y_1),\ldots,1\otimes \rho(y_m))\right ) (q_s\otimes 1) \right \| =0 .\end{align}
 Therefore, by (3.23) and (3.24), there is a positive integer $t$ such that
 \begin{align}
\max_{1\le j\le r} &\left | \| P_j(q_t \pi(x_1) q_t\otimes 1,
\ldots,  q_t \pi(x_n) q_t\otimes 1, q_t\otimes
\rho(y_1),\ldots,q_t\otimes \rho(y_m))\|  \right .\notag\\ &\qquad
\quad - \left . \|P_j(\pi(x_1)\otimes 1, \ldots,\pi(x_n) \otimes 1,
1\otimes\rho( y_1),\ldots,1\otimes \rho(y_m))\|\right |\le \epsilon
/2.\end{align} Moreover, we know
$$
\{q_t \pi(x_1) q_t\otimes 1, \ldots,  q_t \pi(x_n) q_t\otimes 1,
q_t\otimes \rho(y_1),\ldots,q_t\otimes \rho(y_m)\} \subseteq
B(q_t\mathcal H_2)\otimes \rho(\mathcal B).
$$
Because $\mathcal B$ is an MF algebra, $B(q_t\mathcal H_2)\otimes
\rho(\mathcal B)$ is an MF algebra. Thus there are a positive
integer $k$ and
 $$
 \{A_1,\ldots,A_n,B_1,\ldots,B_m\}\subseteq \mathcal M_{k}(\Bbb C)
 $$ such that
 $$\begin{aligned}
  \max_{1\le j\le r} &\left | \|P_j(A_1,\ldots,A_n,B_1,\ldots,B_m)\| \right . \\ & \qquad -
  \|P_j(q_t\pi( x_1) q_t\otimes 1, \ldots,  q_t \pi(x_n) q_t\otimes 1, q_t\otimes\rho( y_1),\ldots,q_t\otimes \rho(y_m))\| \left . \right |\le \epsilon/2,\end{aligned}
 $$ Combining with (3.25), we get
 $$\begin{aligned}
  \max_{1\le j\le r} &\left | \|P_j(A_1,\ldots,A_n,B_1,\ldots,B_m)\| \right .-
  \|P_j(x_1\otimes 1,\ldots,x_n\otimes 1, 1\otimes y_1,\ldots,1\otimes y_m)
  \|_{\mathcal A\otimes_{min}\mathcal B}\left . \right |\le \epsilon.\end{aligned}
 $$ This completes the proof of the proposition.

\end{proof}

\begin{example}Note $C^*(F_n)$ is a residually finite dimensional C$^*$-algebra, thus a quasidiagonal C$^*$-algebra.
 Therefore we have the following result: Suppose $\mathcal B$ is a unital  separable MF algebra. By Corollary 5.3 in \cite{HaSh4},
$C_r^*(F_n)\ast_{\Bbb C}\mathcal B$ is an MF algbebra, where
$C_r^*(F_n)\ast_{\Bbb C}\mathcal B$ is the unital full free
product of $C_r^*(F_n) $ and $ \mathcal B$.
  Then by Proposition 3.1, $C^*(F_n)\otimes_{min}(C_r^*(F_n)\ast_{\Bbb C}\mathcal B)$ is an MF
algebra. Thus, $Ext(C^*(F_n)\otimes_{min}(C_r^*(F_n)\ast_{\Bbb
C}\mathcal B))$ is not a group because of Lemma 2.4.
\end{example}

\section{Crossed Products}

In this section, we   assume that $\mathcal A$ is an MF algebra and
$G$ is a discrete countable amenable group such that there is a
homomorphism $\alpha: G\rightarrow Aut(\mathcal A)$. We will
consider the question when the reduced, or full, crossed product of
$\mathcal A$ by $G$ is an MF algebra. This is not always the case if
we are not going to add extra conditions on the actions of $G$ on
$\mathcal A$. Note Cuntz algebra is stably isometric to the crossed
product of an AF algebra by $\Bbb Z$. Obviously Cuntz algebra is not
an MF algebra but AF algebra is.

Because $G$ is an amenable group, the reduced crossed product of
$\mathcal A$ by $G$ coincides with the full crossed product of
$\mathcal A$ by $G$, which we will denote by $\mathcal
A\rtimes_\alpha G$.
\subsection{Finite group action} When $G$ is a finite group, we have
the following result.
\begin{theorem}
Suppose $\mathcal A$ is an MF algebra and $G$ is a finite group such
that there is a homomorphism $\alpha: G\rightarrow Aut(\mathcal A)$.
Then
 $\mathcal A\rtimes_{\alpha} G$ is an MF algebra.
\end{theorem}

\begin{proof}
Since $G$ is a finite group, by Green's result in \cite{Gr}, we know that
$$
\mathcal A\rtimes_{\alpha} G \hookrightarrow\mathcal
A\otimes_{min} C(G) \rtimes_{\alpha\otimes \gamma} G\simeq
\mathcal A \otimes_{min} \mathcal K,
$$ where $C(G)$ is the C$^*$-algebra consisting all (continuous)
functions on $G$, $\gamma$ is  the action induced by left
translation by $G$ on $C( G) $, and $\mathcal K$ is the set of
compact operators on $l^2(G)$. By the fact the $\mathcal A$ is an
MF algebra and  results in \cite{BK} or Proposition 3.1, $\mathcal
A \otimes \mathcal K$ is an MF algebra. Hence $ \mathcal
A\rtimes_{\alpha} G$ is an MF algebra.
\end{proof}

\begin{corollary}
If $G$ is a
  finite group with $|G|\ge 2$, where $|G|$ is the order of group $G$, then $C^*_r(\Bbb Z \ast G)$ is
an MF algebra and $Ext(C^*_r(\Bbb Z \ast G))$ is not a group.
\end{corollary}
\begin{proof}
Since a  finite group is a subgroup of a finite permutation group,
we need only to prove the result when $G$ is a finite permutation
group $S_n$. If $\alpha: S_n\rightarrow Aut(F_n)$ is a mapping
defined by
$$
\alpha(\sigma)(g_i)=g_{\sigma(i)},
$$ where $g_1,\ldots, g_n$ are   canonical generators of the free
group $F_n$, then $\alpha$ induces an action $S_n$ on $C_r^*(F_n)$
satisfying
$$
C_r^*(F_n)\rtimes_\alpha S_n\simeq C_r^*(F_n\rtimes_\alpha S_n).
$$ Since Haagerup and Thorbj{\o}rnsen proved that $C_r^*(F_n)$ is an MF algebra,
 it follows immediately from last proposition that $C_r^*(F_n)\rtimes_\alpha S_n$
 is an MF algebra.

 Let $g=(g_1g_2\cdots g_n)^3$ be an element in $F_n$. We observe that
  $g$ and $S_n$ are free in $F_n\rtimes_\alpha S_n.$  To see this,
let $\sigma_1,\ldots,\sigma_m\in S_n $ and $n_1,\ldots, n_m\in\Bbb
N$. Then
$$\begin{aligned}
g^{n_1}&\sigma_1g^{n_2} \sigma_2\cdots g^{n_m} \sigma_m \\&=
g^{n_1}(\sigma_1g^{n_2}\sigma_1^{-1})
(\sigma_1\sigma_2)g^{n_3}(\sigma_1\sigma_2)^{-1} \cdots
(\sigma_1\sigma_2\cdots \sigma_{m-1})g^{n_m} (\sigma_1\sigma_2\cdots
\sigma_{m-1})^{-1} (\sigma_1\sigma_2\cdots \sigma_{m-1}\sigma_m)\\
&=g^{n_1} (\beta_1 g^{n_2}\beta_1^{-1})\cdots (\beta_{m-1}g^{n_m}
\beta_{m-1}^{-1}) \beta_m,
\end{aligned}
$$ where
$$
\beta_j=\sigma_1\cdots \sigma_j, \qquad \forall \ 1\le j\le m.
$$Apparently, to check that $g$ and $S_n$ are free in $F_n\rtimes_\alpha
S_n,$ we need only to show that
$$g^{n_1} (\beta_1 g^{n_2}\beta_1^{-1})\cdots (\beta_{m-1}g^{n_m}
\beta_{m-1}^{-1})\ne e, \quad  \forall n_1,\ldots, n_m\in\Bbb
Z\setminus \{0\}, \   \forall \ e\ne\beta_1\ne \beta_2\ne \cdots \ne
\beta_{m-1}\in S_n,$$i.e.
$$
g^{n_1} (\alpha(\beta_1) (g))^{n_2} \cdots (\alpha(\beta_{m-1})
(g))^{n_m}
 \ne e, \quad  \forall n_1,\ldots, n_m\in\Bbb Z\setminus \{0\}, \
\forall \ e\ne\beta_1\ne \beta_2\ne \cdots \ne \beta_{m-1}\in S_n
$$
 Note $g=(g_1g_2\cdots g_n)^3$. It is not
hard to show, using mathematical induction on $m$, that the reduced
form of the word  $g^{n_1} (\alpha(\beta_1) (g))^{n_2} \cdots
(\alpha(\beta_{m-1}) (g))^{n_m}$ in $F_n$ always ends in the letters
of
$$
 g_{_{\beta_{m-1}(1)}}g_{_{\beta_{m-1}(2)}}\cdots g_{_{\beta_{m-1}(n)}}  \qquad
 or \qquad
g_{_{\beta_{m-1}(n)}}^{-1}\cdots g_{_{\beta_{m-1}(2)}}^{-1}
 g_{_{\beta_{m-1}(1)}}^{-1},$$ and it can never be equal to $e$.
Therefore $g$ and $S_n$ are free in $F_n\rtimes_\alpha S_n.$.

  Hence $\Bbb Z\ast
  S_n$ can be viewed as a subgroup of $F_n\rtimes_\alpha S_n$ and $C_r^*(\Bbb Z\ast
  S_n)$ is a C$^*$-subalgebra of $C_r^*(F_n\rtimes_\alpha S_n)$.
  It follows $C_r^*(\Bbb Z\ast
  S_n)$ is also an MF algebra.
  Moreover since $\Bbb Z\ast S_n$ is not an
amenable group for $n\ge 2$, $C_r^*(\Bbb Z\ast S_n)$ is not a
quasidiagonal C$^*$-algebra for $n\ge 2$ by Lemma 2.5, whence
$Ext(C^*_r(\Bbb Z \ast S_n))$ is not a group by Lemma 2.4. So,
for any  finite group $G$ with $|G|\ge 2$, $C^*_r(\Bbb Z \ast G)$
is an MF algebra and $Ext(C^*_r(\Bbb Z \ast G))$ is not a group.

\end{proof}

\begin{remark}
Assume $n\ge 1, m\ge 0$ are nonnegative integers and $r=n+m$. Let
$F_r$ be the free group on $r$ many generators, $g_1,\ldots, g_{n
}, g_{n +1}, \cdots, g_r$. Let $S_n$ be the permutation groups on
the integers $\{1,\ldots,n\}$.

Let $\alpha$ be a homomorphism from $S_n$ into $Aut(F_r)$ defined
the following mapping: $\forall \ \sigma\in S_n$
$$
\begin{aligned}
  &\alpha(\sigma)(g_i)=g_{\sigma(i )} \qquad \quad for \quad    \ \  1 \le i\le
   n\\
  &\alpha(\sigma)(g_j)=g_{j} \qquad \quad for \quad n+1\le j\le
  r.
\end{aligned}
$$Let $F_m$ be the free subgroup generated by $g_{n+1}, \cdots, g_r$ in
$F_r\rtimes_\alpha S_n$. Then $S_n$ and $F_m$, as subgroups of
$F_r\rtimes_\alpha S_n$, commute. Let $h_i= (g_1 ^ig_2 ^i\cdots g_r
^i)^3$, for $1\le i\le n$, be elements in $F_r\rtimes_\alpha S_n$
and $F_n$ be the free subgroup generated by $h_1,\ldots,h_n$ in
$F_r\rtimes_\alpha S_n$. Then, similar to Corollary 4.1, we observe
$F_n$ and $F_m\cup S_n$ are free in $F_r\rtimes_\alpha S_n$.
 Now we have the following observation:
$$
 F_n\ast (S_n\times F_m) \subseteq F_r\rtimes_\alpha S_n.
$$
By Theorem 4.1 and the same arguments as in Corollary 4.1, we have
the following result.
\end{remark}

\begin{corollary}
Assume $n\ge 1, m\ge 0$ are nonnegative integers and $G$ is a
nontrivial finite group. Then $C_r(F_n\ast (G\times F_m))$ is an
MF algebra and $Ext(C_r(F_n\ast (G\times F_m)))$ is not a group
by Lemma 2.4 and Lemma 2.5.
\end{corollary}
\begin{remark}
Assume that $n\ge 2$ is a positive integer and $H_1, \ldots, H_n$ is
a family of finite groups. Let
$$
H= H_1\ast H_2 \ast \cdots \ast H_n
$$ be the free product of the groups  $H_1, \ldots, H_n$. Then  there is a finite group
$G$ satisfying
$$
H \subseteq \Bbb Z \ast G.
$$ Thus
$C_r^*(H)$ is a C$^*$-subalgebra of $C_r^*(\Bbb Z \ast G)$. It
follows from the Corollary 4.1 that $C_r^*(H)$ is an MF algebra.
Hence we have the following result.
\end{remark}
\begin{corollary}
Assume that $n\ge 2$ is a positive integer and $H_1, \ldots, H_n$ is
a family of finite groups. Let
$$
H= H_1\ast H_2 \ast \cdots \ast H_n.
$$then $C_r^*(H)$ is an MF algebra. Moreover, if there are $1\le i\ne j\le n$ such that
$$ |H_i|\ge 2 \qquad and \qquad |H_j|\ge 3,$$ where $|H|$ denotes the order
of the group $H$, then $Ext(C_r^*(H))$ is not a group.
\end{corollary}
\begin{proof}
By the explanation in the Remark 4.2, we know that  $C_r^*(H)$ is an
MF algebra. if there are $1\le i\ne j\le n$ such that $  |H_i|\ge 2
$ and $ |H_j|\ge 3,$  then $H$ is not an amenable group. It follows
that $C_r^*(H)$ is not a quasidiagonal C$^*$-algebra, whence
$Ext(C_r^*(H))$ is not a group.
\end{proof}

\begin{example}
For any positive integers $p\ge 2, q\ge 3$, $Ext(C_r^*(\Bbb Z_p\ast
\Bbb Z_q))$ is not a group.
\end{example}

\subsection{Subgroup of finite index} We have the following result.
\begin{proposition}
Assume that $G$ is a discrete countable group and $H$ is a subgroup
of $G$ with a finite index. If $C_r^*(H)$ is an MF algebra, then
$C_r^*(G)$ is also an MF algebra.
\end{proposition}
\begin{proof} 

Let $l^2(G)$ be
 the Hilbert space with an orthonormal basis $\{e_g\}_{g\in   G}$.
 Recall the left regular representation of the group $G$ is the
 unitary representation $\lambda: G\rightarrow \mathcal U(l^2(G))$
 given by $\lambda(g)(e_{h})=e_{gh}$, $g\in G$. Then $C_r^*(G)$ is the C$^*$-algebra generated by
 $\{\lambda(g)\}_{g\in G}$ in $B(l^2(G))$ and the  C$^*$-subalgebra $\mathcal B$ generated by
 $\{\lambda(g)\}_{g\in H}$ in $C_r^*(G)$ is $*$-isomorphic to $C_r^*(H)$.

 Note $[G:H]<\infty$. There are elements $e=g_1 , g_2,\ldots, g_n$ in
 $G$ such that $G=\cup_{i} g_iH$ and $g_iH\cap g_j H=\emptyset$ if
 $i\ne j$. Let $l^2(H)$ be the closed subspace spanned by $\{e_h\}_{ h\in
 H}$ in  $l^2(G)$ and $E$ be the projection from  $l^2(G)$ onto
 $l^2(H)$. Then for each $1\le i\le n,$
 $\lambda(g_i)   E  \lambda(g_i^{-1})$ is the projection from $l^2(G)$
 onto the closed subspace $l^2(g_iH)$. Moreover
 $$
\sum_{i=1}^n \lambda(g_i)   E  \lambda(g_i^{-1}) =I.
 $$

Note that $ E\mathcal B E\simeq C_r^*(H)  $ and $ \{\lambda(g_i) E
\lambda(g_j^{-1}) \}_{1\le i,j\le n} $ consists a system of
matrix units of $\mathcal M_n(\Bbb C)$. Let $\mathcal A$ be the
C$^*$-algebra generated by elements $ \{\lambda(g_i)   E\mathcal
B E  \lambda(g_j^{-1}) \}_{1\le i,j\le n} $ in  $B(l^2(G))$. Then
we observe that
$$
\mathcal A \simeq C_r^*(H)\otimes \mathcal M_n(\Bbb C).
$$

On the other hand, for any $g\in G$, there is  a permutation
$\sigma_g$ on the integers $\{1,2,\ldots, n\}$ such that
$$
gg_iH=g_{\sigma_{_{g}}(i)} H, \qquad \forall \ 1\le i\le n,\qquad
i.e. \qquad  g_{\sigma_{_{g}}(i)}^{-1}   gg_i \in  H, \qquad
\forall \ 1\le i\le n.
$$
Thus there is some $h_i\in H$ such that $$g_{\sigma_{_{g}}(i)}^{-1}
  gg_i=h_i, \qquad i.e. \qquad
\lambda(g_{\sigma_{_{g}}(i)}^{-1})
 \lambda (g)   \lambda(g_i)=\lambda (h_i). $$ Moreover,
$$  E \lambda(g_{j}^{-1}) \lambda (g) \lambda
(g_i) E=0, \qquad \ if \ \  j\ne \sigma_{_{g}}(i). $$ Therefore
$$
\begin{aligned}
\lambda(g) & = \left ( \sum_{j=1}^n \lambda(g_j)   E
\lambda(g_j^{-1}) \right) \lambda (g) \left ( \sum_{i=1}^n
\lambda(g_i)   E  \lambda(g_i^{-1}) \right) \\
&= \sum_{i,j=1}^n \lambda(g_j)   E  \lambda(g_j^{-1})
 \lambda (g)   \lambda(g_i)   E
\lambda(g_i^{-1})\\
&=   \sum_{i =1}^n \lambda(g_{\sigma_{_{g}}(i)})   E
\lambda(g_{\sigma_{_{g}}(i)}^{-1})
 \lambda (g)   \lambda(g_i)   E
\lambda(g_i^{-1})\\
&= \sum_{i =1}^n \lambda(g_{\sigma_{_{g}}(i)})   E \lambda(h_i)
E  \lambda(g_i^{-1}) \in \mathcal A.
\end{aligned}
$$
It follows that $ C_r^*(G) \subseteq \mathcal A. $ Apparently,
$\mathcal A$ is an MF algebra, because $C_r^*(H)$ is an MF algebra.
Hence, being a subalgebra of an MF algebra,  $C_r^*(G)$ is also an
MF algebra.
\end{proof}

\begin{corollary}
Let $SL_2(\Bbb Z)$ to be the spacial linear group of $2\times 2$
matrices with integer entries. Then $C_r^*(SL_2(\Bbb Z))$ is an MF
algebra and $Ext(C_r^*(SL_2(\Bbb Z)))$ is not a group.
\end{corollary}
\begin{proof}
Note $SL_2(\Bbb Z)$ has a subgroup, which is isomorphic to $F_2$,
with an index $12$. The result of the corollary follows immediately
from the preceding proposition and Haagerup and Thorbj{\o}rnsen's
result that  $C_r^*(F_2)$ is an MF algebra.
\end{proof}



\subsection{Integer group action}

In \cite{PiVoi}, Pimsner and Voiculescu proved the following result:
{\em Suppose $\mathcal A$ is a unital sepaprable C$^*$-algebra and
$\alpha: \Bbb Z\rightarrow Aut(\mathcal A)$ is a homomorphism, such
that there exists a sequence of integers $0\le n_1<n_2<\cdots,$ so
that
$$
\lim_{j\rightarrow \infty} \|\alpha(n_j) a -a\|=0
$$ for all $a\in \mathcal A$. Assume that $\mathcal A$ is quasidiagonal. Then $\mathcal A\rtimes_\alpha \Bbb Z$ is also quasidiagonal.} In this subsection, we will prove an analogue of Pimsner and Voiculescu's result in the context of MF algebras.

Let $l^2(\Bbb Z)$ be
 the Hilbert space with an orthonormal basis $\{e_n\}_{n\in \Bbb Z}$.
 Recall the left regular representation of the group $\Bbb Z$ is the
 unitary representation $\lambda: \Bbb Z\rightarrow \mathcal U(l^2(\Bbb Z))$
 given by $\lambda(n)(e_{n_1})=e_{n+n_1}$, $n_1\in \Bbb Z$. Then $C_r^*(\Bbb Z)$ is the C$^*$-algebra generated by
 $\{\lambda(n)\}_{n\in \Bbb Z}$ in $B(l^2(\Bbb Z))$.
Let $u=\lambda(1)$ be the canonical generator of $C^*_r(\Bbb Z)$, a
bilateral shift operator.

Suppose that $\mathcal A$
 is a unital MF algebra generated by a family of self-adjoint elements
  \begin{align}x_1,\ldots, x_m\in \mathcal A. \end{align}
Suppose
 $\alpha$ is a homomorphism from $\Bbb Z$ to $Aut(\mathcal A)$ such
 that
there exists a sequence of integers $0\le n_1<n_2<\cdots,$
satisfying
$$
\lim_{j\rightarrow \infty} \|\alpha(n_j) a -a\|=0\qquad \text {
for all $a\in \mathcal A$.}
$$  Without loss of generality, we assume that
\begin{align}
\|\alpha(n_j)x_i-x_i\|<\frac 1 j \ , \qquad \forall \ j\ge 1, \
1\le i\le m.
\end{align}
Assume
 $\mathcal A$ acts on a Hilbert space $\mathcal H$.
 Define $\sigma:
 \mathcal A\rightarrow B(\mathcal H\otimes l^2(\Bbb Z))$ to be
 $$
\sigma(a) (\xi\otimes e_n)= (\alpha(-n)a)\xi\otimes e_n, \qquad
\forall \ \ n\in \Bbb Z.
 $$
Then the C$^*$-algebra    generated by $\{\sigma(a)\}_{a\in\mathcal
A}\cup \{I_{\mathcal H}\otimes \lambda(n)\}_{n\in \Bbb Z}$ in
$B(\mathcal H\otimes l^2(\Bbb Z))$ is the
reduced crossed product $\mathcal A\rtimes_{\alpha} \Bbb Z$. 

Following the notation from \cite{PiVoi}, we let for
$j=1,2,\ldots$
\begin{align}f_{k,j}=\cos \frac {k\pi}{2n_j}  \cdot e_{k}+\sin \frac
{k\pi}{2n_j}\cdot e_{k-n_j}\qquad \text{ for $0\le k\le
n_j$}\end{align}
   be a unit vector in $l^2(\Bbb Z)$, and $q_{k,j}$ be the rank one projection from $l^2(\Bbb Z)$ onto  $f_{k,j}$. Let
   $$
      q_j=q_{0,j}+\ldots+ q_{n_j-1,j}.
   $$
   Then by \cite{PiVoi}\begin{align}
   q_j\rightarrow I_{l^2(\Bbb Z)} \ in \ SOT \ \quad and \quad \lim_{j\rightarrow\infty}\|uq_j-q_ju\|_{B(l^2(\Bbb Z))}=0.\end {align}

Let $ y_i=\alpha(-1) x_i\in \mathcal A$, for $1\le i\le m$. We
have the following result.
\begin{lemma} Let $s,t,r,N_1,\ldots, N_r$ be positive integers.
 For any $\epsilon>0$,  $\{G_i\}_{i=1}^s$ in $\Bbb C\langle Z\rangle$,
 $\{H_i\}_{i=1}^t$ in $\Bbb C\langle X_1,\ldots,X_m, Y_1,\ldots, Y_m\rangle$ and
 $\{\{P_j^{(i)}\}_{j=-N_i}^{N_i}\}_{i=1}^r$
in $\Bbb C\langle X_1,\ldots, X_m\rangle$, there are a positive
integer $k$ and $k\times k$ matrices
 $$
 A_1,\ldots,A_m, B_1,\ldots,B_m, U \in \mathcal M_k(\Bbb C)
 $$ such that (i) $$ \|U^*A_i- B_iU^*\| \le \epsilon ,  \ \quad  1\le i\le
 m;$$
 (ii) $$\left | \|G_i(U)\|-\|G_i(u)\|_{C^*_r(\Bbb Z)}\right |\le \epsilon,\quad 1\le i\le s ;$$   (iii)
  $$
 \left |\| H_i(A_1,\ldots,A_m,B_1,\ldots,B_m)\|-\|H_i(x_1,\ldots,x_m,y_1,\ldots,y_m)\|_{\mathcal A}\right | \le \epsilon,\quad 1\le i\le
 t;
  $$
and (iv) $\forall \ \ 1\le i\le r$
$$
\|\sum_{j=-N_i}^{N_i}
P_j^{(i)}(\sigma(x_1),\ldots,\sigma(x_m))(1\otimes
u)^{j}\|-\epsilon\le \|\sum_{j=-N_i}^{0}
P_j^{(i)}(A_1,\ldots,A_m)(U^*)^{-j}+\sum_{j=1}^{N_1}
P_j^{(i)}(A_1,\ldots,A_m) U ^{j}\|.
$$

\end{lemma}

\begin{proof} Without loss of generality, we assume that each $\|x_i\|\le 1$ for $1\le i\le m$.
Let
\begin{align}
x_{k,i}=\alpha(-k)x_i, \ \ y_{k,i}=\alpha(-k)y_i\qquad \ for \ 1\le i\le m, \ k\in \Bbb Z.
\end{align}
By (4.2), we have
$$
\|x_{k+n_j,i}-x_{k,i}\|<\frac 1 j \qquad and \qquad
\|y_{k+n_j,i}-y_{k,i}\|<\frac 1 j \qquad for \ \ 1\le i\le m, \
and \  k,j\ge 1.
$$

For the purpose of simplicity, we will assume that $r=1$. The
general case when $r\ge 1$ follows from the similar argument. Thus
we let
\begin{equation}
P(X_1,\ldots,X_m,Z) = \sum_{i=-N_1}^{0}
P_i^{(1)}(X_1,\ldots,X_m)(Z^*)^{-i}+\sum_{i=1}^{N_1}
P_i^{(1)}(X_1,\ldots,X_m)Z^{i}
\end{equation}
be in $\Bbb C\langle X_1,\ldots,X_m,Z\rangle$  and
$$
M=\sum_{i=-N_1}^{N_1} \| P_i^{(1)}( x_1 ,\ldots, x_m
)\|_{\mathcal A} .
$$

For any $\epsilon>0$, there are a positive integer $L\in \Bbb N$
and a unit vector
\begin{equation}\eta =\sum_{l=-L}^L \xi_{l} \otimes
e_l \end{equation}  in $\mathcal H\otimes l^2(\Bbb Z)$, where each
vector  $\xi_{l} $ is in $\mathcal H$, such that
\begin{equation}
\|P (\sigma(x_1),\ldots,\sigma(x_m),1\otimes u )\|_{B(\mathcal
H\otimes l^2(\Bbb Z))} -2\epsilon\le \|P
(\sigma(x_1),\ldots,\sigma(x_m),1\otimes u)\eta \|.
\end{equation}

For such $\epsilon$ and $\eta$, by (4.3) and (4.4), there is a
positive integer $j>8/ \epsilon$ and a finite rank projection
$$
q_j= q_{0,j}+\ldots+ q_{n_j-1,j}
$$ on $l^2(\Bbb Z)$ such that
\begin{align}
 |\cos\frac {(n_j-1)\pi}{2n_j}|+\sum_{i=-N_1}^{N_1}\sum_{l=-L}^L
(M+1)\|e_{l+i}-\cos\frac{(l+i)\pi}{2n_j}f_{l+i,j}\|
    \le  \frac \epsilon 8;\\
  \|(1\otimes q_j)\eta-\eta\| \le \frac \epsilon {16(M+1)(N_1)^2};\\
 \|uq_j-q_ju\|_{B(l^2(\Bbb Z))} \le \frac\epsilon
{16(M+1)(N_1)^2};
 \\
    \left | \|G_i(q_juq_j)\|_{B(l^2(\Bbb Z))}-\|G_i(u)\|_{B(l^2(\Bbb Z))}\right |
       \le\frac \epsilon 8,  \qquad   1\le i\le s,\end{align}
where $M,N_1,L$ are introduced as above.

 Since $\mathcal A$ is an MF
algebra with a family of elements $ \{x_i, y_i, x_{k,i},
y_{k,i}\}_{1\le i\le m, |k|\le n_j+1} $ in $\mathcal A$, for such
$\epsilon, j$, by Lemma 2.3 there exists a family of
quasidiagonal operators
$$
\{a_i, b_i, a_{k,i}, b_{k,i}\}_{1\le i\le m, |k|\le
n_j+1}\subseteq B(\mathcal H)
$$ satisfying, by (4:b) of Lemma 2.3,
\begin{align}
  \max_{1\le i\le t, |k|\le n_j} &\left | \|H_i(a_{k,1},\ldots,a_{k,m},b_{k,1},\ldots,b_{k,m} )\|-\|H_i(x_{k,1},\ldots,x_{k,m},y_{k,1},\ldots,y_{k,m})\|  \right|\le \frac \epsilon 8\\
 &\qquad \qquad  \max_{1\le i\le m, |k|\le n_j}   \| a_{k+1,i}-b_{k,i}\|   \le \frac \epsilon {
 8}\\
  &\max_{1\le i\le m, |k|\le n_j}\{\|a_{k+n_j,i}-a_{k,i}\|,  \|b_{k+n_j,i}-b_{k,i}\|\}<\frac\epsilon 8\\&\max_{1\le i\le N_1,|k|\le n_j}
 |\|P_i^{(1)}(x_1,\ldots,x_m)\|-\|P_i^{(1)}(a_{k,1},\ldots,a_{k,m})\||\le
 \epsilon/8;
\end{align} and by (4:c) of Lemma 2.3,
\begin{equation}
 \sum_{|k|\le n_j+1} \sum_{-N_1\le i\le N_1} \sum_{  1\le l\le L}\| P_i^{(1)}(a_{k,1},\ldots, a_{k,m})\xi_l-
  P_i^{(1)}(x_{k,1},\ldots, x_{k,m})\xi_l
  \|\le \epsilon /8.
\end{equation}

For such $\epsilon, j$, we let
$$\begin{aligned}
\tilde A_i  &= \sum_{0\le k\le n_j-1} a_{k,j} \otimes q_{k,j},\qquad 1\le i\le m\\
\tilde B_i  &= \sum_{0\le k\le n_j-1} b_{k+1,j} \otimes q_{k,j},\qquad 1\le i\le m\\
\tilde U &=  1\otimes q_j u q_j
\end{aligned}$$
By the choice of   $\eta$ in (4.7) and the fact $u=\lambda(1)$,
we know that
\begin{align}
  P(\sigma(x_1),\ldots,\sigma(x_m),1\otimes u)\eta
  &=\left (\sum_{i=-N_1}^{N_1}P_i^{(1)}(\sigma(x_1),\ldots,\sigma(x_m))(1\otimes
  u)^i \right )\left ( \sum_{l=-L}^L \xi_{l} \otimes
e_l \right )\notag\\
&= \sum_{i=-N_1}^{N_1} \sum_{l=-L}^L \left ( P_i^{(1)}(
\alpha(-i-l)x_1 ,\ldots, \alpha(-i-l)x_m ) \xi_l
\right ) \otimes e_{l+i}\notag\\
&=\sum_{i=-N_1}^{N_1} \sum_{l=-L}^L \left ( P_i^{(1)}(
x_{i+l,1},\ldots, x_{i+l,m} ) \xi_l \right ) \otimes e_{l+i}
\end{align}
Note  $\{q_{k,j}\}_{k=0}^{n_j-1}$ is a family of mutually orthogonal
projections. By inequality (4.16), we know
\begin{align}
\|P_i^{(1)}(\tilde A_1,\ldots, \tilde A_m)\|&=\sup_{0\le k\le
n_j-1}\|P_i^{(1)}(a_{k,1},\ldots, a_{k,m})\|\notag\\&\le \sup_{0\le
k\le n_j-1}\|P_i^{(1)}(x_{k,1},\ldots, x_{k,m})\|+1\le M+1
\end{align}
for all $ -N_1\le i\le N_1. $  By inequality (4.11) and (4.19), we
know that

\begin{align}
 \| &P(\tilde A_1,\ldots, \tilde A_m, \tilde U) -   P(\tilde
A_1,\ldots, \tilde A_m, 1\otimes u)(1\otimes q_j)\|\notag\\
&=\|\left (\sum_{i=-N_1}^{0}
    ( P_i^{(1)}( \tilde A_1,\ldots, \tilde A_m
) (1\otimes
   \tilde U^*)^{- i}+\sum_{i=1}^{N_1}
    ( P_i^{(1)}( \tilde A_1,\ldots, \tilde A_m
) (1\otimes
   \tilde U)^i\right ) \notag\\
   &\ \ \ -\left (\sum_{i=-N_1}^{0}
    ( P_i^{(1)}( \tilde A_1,\ldots, \tilde A_m
)(1\otimes u^*)^{-i}+\sum_{i=1}^{N_1}
    ( P_i^{(1)}( \tilde A_1,\ldots, \tilde A_m
)(1\otimes u)^i\right )(1\otimes q_j) \|\notag
 \\ &\le\sum_{i=-N_1}^{0}
    \|P_i^{(1)}( \tilde A_1,\ldots, \tilde A_m)
\|\|(1\otimes
   q_ju^*q_j)^{-i}-(1\otimes u^*)^{-i}(1\otimes q_j)\| \notag\\&\ \ \quad + \sum_{i=1}^{N_1}
    \| P_i^{(1)}( \tilde A_1,\ldots, \tilde A_m)
\|\|(1\otimes
   q_juq_j)^i-(1\otimes u)^i(1\otimes q_j) \| \quad \le \ \  \frac \epsilon 8.
\end{align}
And, by (4.10) and (4.19),
\begin{align}
\|P(\tilde A_1,\ldots,  \tilde A_m, 1\otimes u)(1\otimes
q_j)\eta&- P(\tilde A_1,\ldots, \tilde A_m, 1\otimes u) \eta\|
\notag \\&\le \|P(\tilde A_1,\ldots, \tilde A_m, 1\otimes
u)\|\|(1\otimes
q_j)\eta-  \eta\|\notag\\
&\le \sum_{i=-N_1}^{N_1} \| P_i^{(1)}( \tilde A_1,\ldots, \tilde
A_m )\|\|(1\otimes q_j)\eta-  \eta\| \notag\\
&\le \frac \epsilon 8.
\end{align}
 Moreover, by the definitions of $\tilde A_1,\ldots, \tilde
A_m $, we have
\begin{align}
  P(\tilde A_1,&\ldots, \tilde A_m, 1\otimes u) \eta
 = \sum_{i=-N_1}^{N_1} \sum_{l=-L}^L    P_i^{(1)}( \tilde
A_1,\ldots, \tilde A_m )(1\otimes
   u)^i  ( \xi_l
  \otimes e_l )\notag\\
 & = \sum_{i=-N_1}^{N_1} \sum_{l=-L}^L     P_i^{(1)}( \tilde A_1,\ldots,
\tilde A_m )(  \xi_l
  \otimes  e_{l+i } )\notag\\
&= \sum_{i=-N_1}^{N_1} \sum_{l=-L}^L    \sum_{0\le k\le n_j-1} \left
( P_i^{(1)}( a_{k,1},\ldots, a_{k,m})\xi_l \otimes q_{k,j}  e_{l+i }
\right )
    \notag\\ &= \sum_{i=-N_1}^{N_1} \sum_{l=-L}^L     ( P_i^{(1)}(
a_{l+i,1},\ldots, a_{l+i,m})   \xi_l)
  \otimes (\cos\frac {(l+i)\pi}{2n_j} f_{l+i,j })
  \notag\\
  &= \sum_{i=-N_1}^{N_1} \sum_{l=-L}^L     ( P_i^{(1)}(
a_{l+i,1},\ldots, a_{l+i,m})   \xi_l)
  \otimes e_{l+i}   \notag\\
  & \ \ \ \ -  \sum_{i=-N_1}^{N_1} \sum_{l=-L}^L      ( P_i^{(1)}(
a_{l+i,1},\ldots, a_{l+i,m})   \xi_l)
  \otimes (e_{l+i}- \cos\frac {(l+i)\pi}{2n_j} f_{l+i,j } )
\end{align}
By
 inequalities    (4.20), (4.21),  (4.18),   (4.22), (4.17) and (4.9)  we
know,

\begin{align}
\|P(\tilde A_1,\ldots,& \tilde A_m, \tilde U)\eta  - P
(\sigma(x_1),\ldots,\sigma(x_m),1\otimes u)\eta\|\notag\\
\le &\|P(\tilde A_1,\ldots, \tilde A_m, \tilde U)\eta  - P
(\tilde A_1,\ldots, \tilde A_m,1\otimes u)(1\otimes
q_j)\eta\|\notag \\
 & + \|P(\tilde A_1,\ldots, \tilde A_m, 1\otimes u)(1\otimes
q_j)\eta  - P (\tilde A_1,\ldots, \tilde A_m,1\otimes u)
 \eta\|\notag \\
 &  + \|P (\tilde A_1,\ldots, \tilde A_m,1\otimes u) \eta  -P
(\sigma(x_1),\ldots,\sigma(x_m),1\otimes u) \eta\|\notag
 \\ \le &  \frac {\epsilon} 8+\frac \epsilon 8+ \sum_{i=-N_1}^{N_1} \sum_{l=-L}^L
    \|P_i^{(1)}(
a_{l+i,1},\ldots, a_{l+i,m})\xi_l -P_i^{(1)}( x_{l+i,1},\ldots,
x_{l+i,m})\xi_l \|\notag\\
  &\quad \quad +
  \sum_{i=-N_1}^{N_1} \sum_{l=-L}^L      \| P_i^{(1)}(
a_{l+i,1},\ldots, a_{l+i,m}) \| \| e_{l+i}- \cos\frac
{(l+i)\pi}{2n_j} f_{l+i,j } \|\notag \\
 \le &\frac \epsilon 2 . \end{align} Note that   $ \{a_i, b_i,
a_{k,i}, b_{k,i}\}_{1\le i\le m, |k|\le n_j+1}\subseteq
B(\mathcal H) $ is a family of quasidiagonal operators.    There
is a finite rank projection  $p $ on $\mathcal H$ such that, by
Lemma 2.1 part (ii),
\begin{align}
& \max_{1\le i\le t}\left |
\|H_i(a_{k,1},\ldots,a_{k,m},b_{k,1},\ldots,b_{k,m} )\|
 -
  \|H_i(p a_{k,1}p ,\ldots,p a_mp ,p b_{k,1}p ,\ldots, b_{k,m} )
  \|\right |\le \frac \epsilon 8;
\end{align} and, by Lemma 2.1 part (i),
\begin{equation}
\|P( A_1,\ldots,  A_m,  U)\eta- P(\tilde A_1,\ldots, \tilde A_m,
\tilde U)\eta\|\le \frac {\epsilon} 2, \end{equation} where
$$\begin{aligned}
A_i  &=p\tilde A_ip= \sum_{0\le k\le n_j-1}p a_{k,j}p \otimes q_{k,j},\qquad 1\le i\le m\\
B_i  &= p\tilde B_ip=\sum_{0\le k\le n_j-1}p b_{k+1,j}p \otimes q_{k,j},\qquad 1\le i\le m\\
U &= p\tilde Up= p \otimes q_j u q_j
\end{aligned}
$$ are the elements in $B((p\otimes q_j)(\mathcal H\otimes l^2(\Bbb
Z)))$. By (4.8), (4.23) and (4.25), we have
\begin{equation}
\|P( A_1,\ldots,  A_m,  U)\|\ge \| P
(\sigma(x_1),\ldots,\sigma(x_m),1\otimes u) \|_{B(\mathcal
H\otimes l^2(\Bbb Z) )}-4\epsilon.
\end{equation}

Since $\{q_{k,j}\}_{k=0}^{n_j-1}$ is a family of mutually orthogonal
projections, we have
\begin{align}
\|H_i(A_1,\ldots,A_m,B_1,\ldots, B_m)\| =\max_{1\le k\le n_j} \{\|H_i(pa_{k,1}p,\ldots,pa_{k,m}p,pb_{k,1}p,\ldots,pb_{k,m}p)\|\}.
\end{align}
Now it follows from  (4.13),   (4.24) and (4.27) that for $1\le
i\le t$
\begin{align}
\left | \|H_i(A_1,\ldots,A_m,B_1,\ldots, B_m)\|- \|H_i(x_1,\ldots,x_m,y_1,\ldots,y_m)\|_{\mathcal A} \right |\le \epsilon.
\end{align}
By (4.12), we know that for $1\le i\le  s$
\begin{align}
\left |\|G_i(U)\|- \|G_i(u)\|\right|\le \epsilon.
\end{align}

Let $E_n$ be the rank one projection from $l^2(\Bbb Z)$ onto the
vector $e_n$. For any unit vector $h\in \mathcal H$
$$
 \begin{aligned}
  &\left ( a_{k,i}\otimes E_k+a_{k-n_j,i}\otimes E_{k-n_j}\right) h\otimes f_{k,j} \\
  & \qquad = \left( a_{k,i}\otimes E_k+a_{k-n_j,i}\otimes E_{k-n_j} \right)(h\otimes (\cos \frac {k\pi}{2n_j} \cdot e_{k}+
  \sin \frac {k\pi}{2n_j}\cdot e_{k-n_j}))\\
  &\qquad =(\cos\frac {k\pi}{2n_j} \cdot a_{k,i} h )\otimes e_k + (\sin \frac {k\pi}{2n_j} \cdot   a_{k-n_j,i})\otimes e_{k-n_j}\\
  &\qquad = (\cos\frac {k\pi}{2n_j}  \cdot  a_{k,i} h) \otimes e_k+ (\sin \frac {k\pi}{2n_j}  \cdot  a_{k,i})\otimes
  e_{k-n_j}
  +(\sin \frac {k\pi}{2n_j}  \cdot  (a_{k-n_j,i}-a_{k,i}))\otimes e_{k-n_j}\\
  &\qquad =(a_{k,i}\otimes q_{k,j})(h\otimes f_{k,j})+  (\sin \frac {k\pi}{2n_j}  \cdot  (a_{k-n_j,i}-a_{k,i}))\otimes
  e_{k-n_j}; \quad \ \ \forall \ 1\le i\le m.
 \end{aligned}
$$By (4.15), we have that
$$\begin{aligned}
\left \| \left ( a_{k,i}\otimes E_k+a_{k-n_j,i}\otimes
E_{k-n_j}\right)(1\otimes q_{k,j}) - \left (a_{k,i}\otimes
q_{k,j}\right) \right \|\le \epsilon /8; \quad \ \ \forall \ 1\le
i\le m.\end{aligned}
$$
Together with the fact that $\{q_{k,j}\}_{k=0}^{n_j-1}$ is a family
of mutually orthogonal projections, we have
 \begin{align}
&\left \|\left (\sum_{0\le k\le n_k-1} \left ( a_{k,i}\otimes
E_k+a_{k-n_j,i}\otimes E_{k-n_j}\right)\right)(1\otimes q_{ j}) -
\sum_{0\le k\le
n_k-1} \left (a_{k,i}\otimes q_{k,j}\right) \right \| \notag\\
& \qquad \qquad=\left \|\sum_{0\le k\le n_k-1} \left (
a_{k,i}\otimes E_k+a_{k-n_j,i}\otimes E_{k-n_j}\right)(1\otimes
q_{k,j}) - \sum_{0\le k\le
n_k-1} \left (a_{k,i}\otimes q_{k,j}\right) \right \| \notag\\
&\qquad \qquad=\left \|\sum_{0\le k\le n_k-1} \left (\left (
a_{k,i}\otimes E_k+a_{k-n_j,i}\otimes E_{k-n_j}\right)(1\otimes
q_{k,j}) -
  \left (a_{k,i}\otimes q_{k,j}\right) \right )\right \| \notag\\
  &\qquad \qquad\le \max_{0\le k\le n_k-1}\left \|   \left (
a_{k,i}\otimes E_k+a_{k-n_j,i}\otimes E_{k-n_j}\right)(1\otimes
q_{k,j}) -
  \left (a_{k,i}\otimes q_{k,j}\right) \right \|\notag \\
 &\qquad \qquad \le \epsilon/ 8; \quad \ \ \forall \ 1\le i\le m. \end{align}
Similarly,
 \begin{align}
&\left \|\left (\sum_{0\le k\le n_k-1} \left ( b_{k,i}\otimes
E_k+b_{k-n_j,i}\otimes E_{k-n_j}\right)\right)(1\otimes q_{ j}) -
\sum_{0\le k\le n_k-1} \left (b_{k,i}\otimes q_{k,j}\right)
\right \|  \le \epsilon/ 8. \end{align}

On the other hand, we observe that
\begin{align}
  (1\otimes u^*) \sum_{0\le k\le n_j-1} \left ( a_{k,i}\otimes E_k+a_{k-n_j,i}\otimes E_{k-n_j}\right)
  =\sum_{0\le k\le n_j-1} \left ( a_{k,i}\otimes E_{k-1}+a_{k-n_j,i}\otimes E_{k-n_j-1}\right) (1\otimes u^* )
\end{align}
and by (4.9),  for $1\le i\le m,$  \begin{align}
(a_{-n_j,i}\otimes E_{-n_j-1}) (1\otimes q_j)=0; \qquad \|
(b_{n_j-1,i}\otimes E_{n_j-1}) (1\otimes q_j) \| \le |\cos\frac
{(n_j-1)\pi}{2n_j}|\le \epsilon/8.
\end{align}
It follows that for $1\le i\le m$
\begin{align}
T_i=&\left . (1\otimes u^*) \sum_{0\le k\le n_j-1}  \left (
a_{k,i}\otimes
E_k+a_{k-n_j,i}\otimes E_{k-n_j}\right)(1\otimes q_j) \right .\notag\\
&\qquad \qquad \qquad \left. - \left (\sum_{0\le k\le n_j-1}
\left ( b_{k,i}\otimes E_k+b_{k-n_j,i}\otimes
E_{k-n_j}\right)\right
)(1\otimes q_j) (1\otimes u^*)\right .\notag \\
&=\left .  \sum_{0\le k\le n_j-1}  \left ( a_{k,i}\otimes
E_{k-1}+a_{k-n_j,i}\otimes E_{k-1-n_j}\right) (1\otimes u^*)(1\otimes q_j) \right .\tag{By (4.32)}\\
&\qquad \qquad \qquad \left. - \left (\sum_{0\le k\le n_j-1} \left (
b_{k,i}\otimes E_k+b_{k-n_j,i}\otimes E_{k-n_j}\right)\right
)(1\otimes q_j  u^*) \right .\notag
\end{align}
\begin{align}&=\left . \sum_{0\le k\le n_j-1}  \left (
a_{k,i}\otimes
E_{k-1}+a_{k-n_j,i}\otimes E_{k-1-n_j}\right) ((1\otimes q_j  u^*)-1\otimes (q_ju^*-u^*q_j))\right .\notag\\
&\qquad \qquad \qquad \left. - \left (\sum_{0\le k\le n_j-1}
\left ( b_{k,i}\otimes E_k+b_{k-n_j,i}\otimes
E_{k-n_j}\right)\right
)(1\otimes q_j  u^*) \right .\notag  \\
&=\left . - \left (\sum_{0\le k\le n_j-1}  \left ( a_{k,i}\otimes
E_{k-1}+a_{k-n_j,i}\otimes E_{k-1-n_j}\right) \right ) ( 1\otimes (q_ju^*-u^*q_j))\right .\notag\\
&\qquad \qquad \qquad \left. + \left (\sum_{-n_j \le k\le n_j-2}
(a_{k+1,i}-b_{k,i})\otimes E_k   \right )(1\otimes q_j  u^*) \right
.\notag\\& \qquad \qquad \qquad  \left .+ \left( a_{-n_j,i}\otimes
E_{-n_j-1}-b_{n_j-1,i}\otimes E_{n_j-1} \right )(1\otimes q_j
u^*)\right .\notag
\end{align}
By (4.11), (4.14) and (4.33), we have that
\begin{align}
  \|T_i\|\le \frac \epsilon 2.
\end{align}
 Moreover,
\begin{align}
U^*A_i &= (p\otimes  q_ju^* q_j) (\sum_{0\le k\le n_k-1}pa_{k,j}p\otimes q_{k,j})\notag\\
&=(p\otimes q_j )(1\otimes u^*)\left ( \sum_{0\le k\le n_k-1} \left
(a_{k,i}\otimes q_{k,j}\right)
- \sum_{0\le k\le n_k-1} \left ( a_{k,i}\otimes E_k+a_{k-n_j,i}\otimes E_{k-n_j}\right) \right. \notag\\
&\qquad \qquad \qquad \quad\quad \quad \ \qquad +\left .\sum_{0\le
k\le n_k-1} \left ( a_{k,i}\otimes E_k+a_{k-n_j,i}\otimes
E_{k-n_j}\right)\right )(p\otimes q_j)\notag
\end{align} And
\begin{align}
B_iU^* &=   (\sum_{0\le k\le n_k-1}pb_{k,j}p\otimes q_{k,j}) (p\otimes q_j u^* q_j )\notag\\
&=(p\otimes q_j ) \left ( \sum_{0\le k\le n_k-1} \left
(b_{k,i}\otimes q_{k,j}\right) - \sum_{0\le k\le n_k-1} \left (
b_{k,i}\otimes E_k+b_{k-n_j}\otimes E_{k-n_j}\right)  \right
.\notag\\&\qquad \qquad \qquad \qquad +\left .\sum_{0\le k\le n_k-1}
\left ( b_{k,i}\otimes E_k+b_{k-n_j}\otimes E_{k-n_j}\right)\right
)(1\otimes q_j)(1\times u^*)(p\otimes q_j)\notag
\end{align}
By (4.30), (4.31), (4.34), we have
\begin{align}
\|U^*A_i-B_iU^*\|\le    \epsilon  . \end{align} Let
$k=rank(p\otimes q_j)$. Then $B((p\otimes q_j)(\mathcal H\otimes
l^2(\Bbb Z)))\simeq \mathcal M_{k}(\Bbb C)$. By (4.35), (4.29),
(4.28) and (4.26), we have completed the  proof of the lemma.
\end{proof}

Recall $x_1,\ldots, x_m$ is a family of generators of a unital MF
algebra $\mathcal A$; $\alpha$ is an action of $\Bbb Z$ on $\mathcal
A$; $u=\lambda(1)$ in $C_r^*(\Bbb Z)$  is the bilateral shift; and
$y_i=\alpha(-1)x_i$ for $i=1,\ldots, m$. From now on, we let
$\{H_r\}_{r=1}^\infty$ (respectively, $\{G_r\}_{r=1}^\infty$) be the
collection of all polynomials in $\Bbb C\langle X_1,\ldots,X_n,Y_1,
\ldots, Y_m\rangle$ (or $\Bbb C\langle Z\rangle$ respectively) with
rational coefficients.

\begin{lemma}
For any $r_0\in\Bbb N$ and $P_1,\ldots,P_{r_0}$ in $\Bbb C\langle
X_1,\ldots,X_n,Z\rangle$, there is some $r_1>0$ such that the
following is true: for   any $k\in \Bbb N$ and any
            $$
 A_1,\ldots,A_m, B_1,\ldots,B_m, U \in \mathcal M_k(\Bbb C)
 $$ satisfying (i) $$ \|U^*A_i-B_iU^*\|_{\mathcal M_k(\Bbb C)} \le \frac 1 {r_1} ,  \ \quad  1\le i\le m;$$  and (ii) $$\left | \|G_i(U)\|_{\mathcal M_k(\Bbb C)}-\|G_i(u)\|_{C^*_r(\Bbb Z)}\right |\le \frac 1 {r_1},\quad 1\le i\le {r_1} ;$$ and (iii)
  $$
 \left |\| H_i(A_1,\ldots,A_m,B_1,\ldots,B_m)\|_{\mathcal M_k(\Bbb C)}-\|H_i(x_1,\ldots,x_m,y_1,\ldots,y_m)\|_{\mathcal A}\right | \le \frac 1 {r_1},\quad 1\le i\le r_1,
  $$    we have
  $$
   \|P_j(A_1,\ldots,A_m, U)\|_{\mathcal M_k(\Bbb C)}-\|P_j(x_1,\ldots,x_m,u)\|_{\mathcal A\rtimes_{\alpha} \Bbb Z}  \le \frac 1 {r_0}, \quad 1\le j\le r_0.
  $$
\end{lemma}

\begin{proof}
Assume the result of the lemma is not true. Thus there are some $r_0>0$, a sequence of positive integers $\{ k_l\}_{l=1}^\infty$, and matrices
 $$
 A_1^{(l)},\ldots,A_m^{(l)}, B_1^{(l)},\ldots,B_m^{(l)}, U^{(l)} \in \mathcal M_{k_l}(\Bbb C)
 $$ satisfying  \begin{enumerate}\item [(a)] $$ \|(U^{(l)})^*A_i^{(l)}-B_i^{(l)}(U^{(l)})^*\|_{\mathcal M_{k_l}(\Bbb C)} \le \frac 1 { l} ,  \ \quad  1\le i\le m;$$  \item [(b)] $$\left | \|G_i(U^{(l)})\|_{\mathcal M_{k_l}(\Bbb C)}-\|G_i(u)\|_{C^*_r(\Bbb Z)}\right |\le \frac 1 {l},\quad 1\le i\le  l ;$$ \item [(c)]
  $$
 \left |\| H_i(A_1^{(l)},\ldots,A_m^{(l)},B_1^{(l)},\ldots,B_m^{(l)})\|_{\mathcal M_{k_l}(\Bbb C)}-\|H_i(x_1,\ldots,x_m,y_1,\ldots,y_m)\|_{\mathcal A}\right | \le \frac 1 { l},\quad 1\le i\le  l.
  $$  \item [(d)]  but
  $$
  \max_{1\le j\le r_0}\{  \|P_j(A_1,\ldots,A_m, U)\|_{\mathcal M_{k_l}(\Bbb C)}-\|P_j(x_1,\ldots,x_m,u)\|_{\mathcal A\rtimes_{\alpha} \Bbb Z} \} >\frac 1 {r_0}.
  $$\end{enumerate}

Consider the unital C$^*$-algebra $\prod_{l}\mathcal M_{k_l}(\Bbb C)/\sum_{l}\mathcal M_{k_l}(\Bbb C)$. Let
$$
a_i=[\langle A_i^{(l)}\rangle_{l=1}^\infty], \ \ b_i=[\langle B_i^{(l)}\rangle_{l=1}^\infty], \ \ w=[\langle U^{(l)}\rangle_{l=1}^\infty], \ \ 1\le i\le m
$$ be the elements in $\prod_{l}\mathcal M_{k_l}(\Bbb C)/\sum_{l}\mathcal M_{k_l}(\Bbb C)$ and
$C^*(a_1,\ldots,a_m, b_1,\ldots,b_m,w)$ be the C$^*$-algebra generated by
$a_1,\ldots,a_m, b_1,\ldots,b_m,w$ in
$\prod_{l}\mathcal M_{k_l}(\Bbb C)/\sum_{l}\mathcal M_{k_l}(\Bbb C)$. By (b) and (c),
 there exist   embedding
$$
\rho_1: \mathcal A\rightarrow  C^*(a_1,\ldots,a_m, b_1,\ldots,b_m,w)\ \ \ and \ \ \ \rho_2: C^*_r(\Bbb Z)\rightarrow  C^*(a_1,\ldots,a_m, b_1,\ldots,b_m,w)
$$ such that
$$
\rho_1(x_i)=a_i, \ \rho_1(y_i)=b_i,\ \ \ \rho_2(u)=w.
$$ By (a), we know that
$$
w^*a_iw=b_i.
$$ Note that $\mathcal A\rtimes_{\alpha}\Bbb Z$ is also the full crossed product of $\mathcal A$ by $\Bbb Z$. From the universal property of the full crossed product, it follows that there is a $*$-homomorphism $$ \rho: \mathcal A\rtimes_{\alpha}\Bbb Z \rightarrow   C^*(a_1,\ldots,a_m, b_1,\ldots,b_m,w)$$ such that $$
\rho(x_i)=a_i, \ \rho (y_i)=b_i,\ \ \ \rho (u)=w.
$$ It follows that for all $1\le j\le r_0$
$$\begin{aligned}
\|P_j(x_1,\ldots,x_m,u)\|_{\mathcal A\rtimes_{\alpha} \Bbb Z}&\ge
\|P_j(a_1,\ldots,a_m,w)\|_{\prod_{l}\mathcal M_{k_l}(\Bbb
C)/\sum_{l}\mathcal M_{k_l}(\Bbb C)}\\&=\limsup_{l}
\|P_j(A_1,\ldots,A_m, U)\|_{\mathcal M_{k_l}(\Bbb
C)}.\end{aligned}
$$
 This contradicts with the assumption (d). Thus the
proof of the lemma is completed. \end{proof}

Combining Lemma 4.1 and Lemma 4.2, we   have the following
conclusion.

\begin{theorem}
Suppose that $\mathcal A$ is a finitely generated unital MF algebra and $\alpha$ is a homomorphism from $\Bbb Z$ into $Aut(\mathcal A)$ such that there is a sequence of integers $0\le n_1<n_2<\cdots$ satisfying
$$
\lim_{j\rightarrow \infty}\|\alpha(n_j) a-a\|=0
$$ for any $a\in \mathcal A$. Then $\mathcal A\rtimes_{\alpha}\Bbb Z $ is an MF algebra.

\end{theorem}
 \begin{proof}
Assume that $\mathcal A$ is generated by a family of self-adjoint
elements $x_1,\ldots, x_m$ and $u$ is the canonical generator of
$C_r(\Bbb Z)$.  Without loss of generality, we assume that
$x_1,\ldots,x_m,u$ generate $\mathcal A\rtimes_{\alpha}\Bbb Z $. Let
$y_i=\alpha(-1)x_i$ for $i=1,\ldots, m$. Let $r_0$ be a positive
integer.

Assume that $\epsilon<1/r_0$ is a positive number. Assume that
$s,t, N_1,\ldots, N_{r_0}$ are positive integers and
 $\{G_i\}_{i=1}^s\subseteq\Bbb C\langle Z\rangle$,
 $\{H_i\}_{i=1}^t\subseteq\Bbb C\langle X_1,\ldots,X_m, Y_1,\ldots, Y_m\rangle$ and
 $\{\{P_j^{(i)}\}_{j=-N_i}^{N_i}\}_{i=1}^{r_0}\subseteq\Bbb C\langle X_1,\ldots, X_m\rangle$.

For $1\le i\le r_0$, we let
\begin{equation*}
P_i(X_1,\ldots,X_m,Z) = \sum_{j=-N_1}^{0}
P_j^{(i)}(X_1,\ldots,X_m)(Z^*)^{-j}+\sum_{i=1}^{N_1}
P_j^{(i)}(X_1,\ldots,X_m)Z^{j}
\end{equation*}
 By Lemma 4.1, there
are a positive integer $k$ and $k\times k$ matrices
 $$
 A_1,\ldots,A_m, B_1,\ldots,B_m,U  \in \mathcal M_k(\Bbb C)
 $$ such that (i) \begin{equation}\|U^*A_i- B_iU^*\| \le \epsilon \le 1/r_0,  \ \quad  1\le i\le
 m;\end{equation}
 (ii) $$\left | \|G_i(U)\|-\|G_i(u)\|_{C^*_r(\Bbb Z)}\right |\le \epsilon,\quad 1\le i\le s ;$$
    (iii)
  $$
 \left |\| H_i(A_1,\ldots,A_m,B_1,\ldots,B_m)\|-\|H_i(x_1,\ldots,x_m,y_1,\ldots,y_m)\|_{\mathcal A}\right | \le \epsilon,\quad 1\le i\le
 t;
  $$
(iv) $\forall \ \ 1\le i\le r_0$
$$
\|  P_ i (x_1,\ldots, x_m, u)\|_{ \mathcal A\rtimes_{\alpha}\Bbb
Z}-\epsilon\le \|  P_ i (A_1,\ldots,A_m, U)\|.
$$

By Lemma 4.2, we know when $\epsilon$ is small enough and $s,t$
are large enough,
$$
 \| P_ i (A_1,\ldots,A_m, U)\|\le \| P_ i (x_1,\ldots, x_m, u)\|_{ \mathcal A\rtimes_{\alpha}\Bbb
Z}+ 1/r_0, \ \ \forall \ \ 1\le i\le r_0 .
$$ Combining with (iv), we have
\begin{equation}
| \|P_ i (A_1,\ldots,A_m, U) \|- \| P_ i (x_1,\ldots, x_m, u)\|_{
\mathcal A\rtimes_{\alpha}\Bbb Z}|\le 1/r_0, \ \ \forall \ \ 1\le
i\le r_0 .\end{equation} By (4.36), (4.37) and Lemma 2.3, we know
that $ \mathcal A\rtimes_{\alpha}\Bbb Z $ is an MF algebra
 \end{proof}

The following corollary follows directly from the preceding theorem.
\begin{corollary} Let $C_r^*(F_2)$ be the reduced C$^*$-algebra of
the free group $F_2$.  Let $u_1,u_2$ be the canonical unitary
generators of $C_r^*(F_2)$ and $0<\theta  <1$ be a positive
number. Let $\alpha$ be a homomorphism from $\Bbb Z$ into
$Aut(C_r^*(F_2))$ induced by the following mapping: $\forall \
n\in \Bbb Z$,
$$
\alpha(n)(u_1)=e^{2n \pi    \theta \cdot  i} u_1 \qquad and \qquad
\alpha(n)(u_2)=e^{2n\pi   \theta \cdot i} u_2.
$$ Then $C_r^*(F_2)\rtimes_{\alpha}\Bbb Z $ is an MF algebra and
$Ext(C_r^*(F_2)\rtimes_{\alpha}\Bbb Z)$ is not a group.
\end{corollary}

\end{document}